\documentclass[11pt, letterpaper, oneside]{article}

\usepackage{latexsym, amsmath, amssymb, amsfonts, amscd}
\usepackage{amsthm}
\usepackage{t1enc}
\usepackage[mathscr]{eucal}
\usepackage{indentfirst}
\usepackage{graphicx}   
\usepackage{fancyhdr}
\usepackage{fancybox}
\usepackage{enumerate,psfrag}
\usepackage[all,poly,web,knot]{xy}
\usepackage{epsfig, diagrams}

\theoremstyle{plain}
\newtheorem{thm}{Theorem}[section]

\newtheorem{prop}[thm]{Proposition}
\newtheorem{lemma}[thm]{Lemma}

\renewcommand{\latticebody}{\drop@{ }}

\theoremstyle{definition}
\newtheorem{defi}[thm]{Definition}

\theoremstyle{remark}
\newtheorem{remark}[thm]{Remark}
\newtheorem{ep}[thm]{Example}

\newcommand{\la}{\leftarrow}
\newcommand{\ra}{\rightarrow}

\newcommand{\rra}{\Rightarrow}

\newcommand{\thra}{\twoheadrightarrow}
\newcommand{\me}{\stackrel{\sim}{\leftrightarrow}} 
\newcommand{\lhc}{\stackrel{\sim}{\to}} 

\newcommand{\N}{\ensuremath{\mathbb N}}
\newcommand{\Z}{\ensuremath{\mathbb Z}}
\newcommand{\C}{\ensuremath{\mathbb C}}
\newcommand{\R}{\ensuremath{\mathbb R}}

\newcommand{\g}{\ensuremath{\frak{g}}}

\newcommand{\D}{\mathscr D}

\newcommand{\cC}{\mathcal{C}}            
\newcommand{\cF}{\mathcal{F}}
\newcommand{\cX}{\mathcal{X}}
\newcommand{\cY}{\mathcal{Y}}

\newcommand{\cG}{\mathcal{G}}
\newcommand{\cH}{\mathcal{H}}

\newcommand{\cO}{\mathcal{O}}

 \DeclareMathOperator{\im}{Im}

\DeclareMathOperator{\colim}{colim}

\newcommand{\tS}{\tilde{S}}

\newcommand{\tphi}{\tilde{\phi}}

\newcommand{\bt}{\mathbf{t}}                  
\newcommand{\bs}{\mathbf{s}}                  

\newcommand\gm[2]{\xymatrix@C=.45cm{#1 \ar@{.>}[r] & #2}} 

\def\L{\Lambda}
\def\D{\Delta}
\newarrow{into} C--->
\newarrow{Eq} =====
\newarrow{dashto} ....>

\def\pD{\partial\D}

\def\PB(#1,#2,#3,#4){
\left\{\begin{matrix}#1&\!\!\!\stackrel{?}{\longrightarrow}&\!\!\!#2\\
\downarrow&&\!\!\!\downarrow\\
#3&\!\!\!\stackrel{?}{\longrightarrow}&\!\!\!#4\end{matrix}\right\}}

\def\pb(#1,#2,#3,#4){ \hom(#1 \to #3, #2 \to #4)}

\begin{document}

\title{Lie II theorem for Lie algebroids via higher groupoids}
\author{Chenchang Zhu  }
\date{\today}
\maketitle

\begin{abstract}
  Following Sullivan's
  spacial realization of a differential algebra, we construct a
  universal integrating  Lie 2-groupoid for every Lie algebroid. Then We show that unlike Lie algebras which one-to-one correspond to
  simply connected Lie groups, Lie algebroids (integrable or not)
  one-to-one correspond to a sort of \'etale Lie 2-groupoids
  with 2-connected source fibres.  Finally, 
  we discuss how to lift Lie algebroid morphisms to Lie 2-groupoid
  morphisms (Lie II Theorem). 
\end{abstract}

\tableofcontents

\section{Introduction}
Lie II theorem for Lie algebras is about how to lift morphisms of Lie
algebras to morphisms of Lie groups.  A
Lie algebroid is the infinitesimal data for a Lie groupoid, as a
Lie algebra is for a Lie group. More precisely,  a Lie
algebroid over a manifold $M$ is a vector bundle $\pi: A\to M$
  with a real Lie bracket $[\,,\,]$ on its space of
sections $H^0(M,A)$ and a bundle map $\rho :A\to TM$ such that the
Leibniz rule
$$[X,f\cdot Y](x)=f(x)\cdot [X,Y](x)+\rho(X)(f)(x)\cdot Y(x)$$ holds for all $X,Y\in
H^0(M,A)$, $f\in C^\infty(M)$ and $x\in M$. Hence when $M$ is a
point, a Lie algebroid becomes a Lie algebra. Also a tangent bundle
$TM \to M$ is certainly a Lie algebroid with $[\,,\,]$ the Lie bracket
of vector fields. The next example is a Poisson manifold $P$ with
$A=T^*P\to P$ and $[df, dg] = d\{f, g\}$ determined by the Poisson
bracket $\{,\}$ (see the book \cite{cw} for a friendly
introduction). Other fundamental examples come from bundles of Lie
algebras of Douady-Lazard (the case when $\rho=0$); Lie
algebra actions on manifolds by Palais \cite{palais}; abstract Atiyah
sequences of Almeida, Molino \cite{am2} and Mackenzie \cite{Mac}. 

Thus Lie II theorem for Lie algebroids is about how to lift an infinitesimal
morphism  from the level of  Lie
algebroids to a global morphism. Its version having Lie groupoids
as the global objects integrating Lie algebroids is well-known \cite{maxu}
\cite{mm:lie2}. The subject of Lie groupoids itself was widely studied
by many people---Connes, Haefliger, Hilsum,
Pradines \cite{pradines}, Skandalis from 60's to 80's---for
applications to
foliations on manifolds. However,
unlike (finite dimensional) Lie algebras which always have their
associated Lie groups, Lie algebroids do not always have their
associated Lie groupoids \cite{am1} \cite{am2}. The complete
integrability criteria is given in a remarkable work of Crainic and
Fernandes \cite{cf}.  However, we claim the situation is not totally hopeless: if we are willing to
enter the world of stacks, we do have the full  one-to-one
correspondence (Lie III) parallel to the classical one for Lie algebras \cite{tz},
\[ \xymatrix{
& \fbox{\parbox{.3\linewidth}{\center{Lie algebras}}}
  \ar[rrrr]^{\text{differentiation at identity}} &  & & &
  \fbox{\parbox{.3\linewidth}{\center{Lie groups}}}
  \ar[llll]^{\text{integration}}   \\
& \fbox{\parbox{.3\linewidth}{\center{Lie algebroids}}}
  \ar[rrrr]^{\text{differentiation at identity}} &  & & &
  \fbox{\parbox{.3\linewidth}{\center{``{\'etale stacky Lie groupoids}''}}}
  \ar[llll]^{\text{integration}}}
\]
Here an {\em \'etale stacky Lie groupoid} $\cG \rra M$  (we
call this sort of groupoid {\em W-groupoid}, as its existence was
first conjectured by Alan
Weinstein \cite{W3} \cite{picard}) is a groupoid in the category of
differentiable stacks with $\cG$ an \'etale stack and $M$ a manifold
(see \cite{z:tgpd-2} for the exact definition), and we call the
procedure of passing from an infinitesimal object to a global object
``integration''. This problem is already very interesting, as shown by
Cattaneo and Felder \cite{cafe}, in the case of Poisson manifolds: the
object integrating $A=T^*P$ is the phase space of the Poisson sigma model
and is further a {\em symplectic}  W-groupoid \cite{tz2}.  

Therefore my effort in this paper is to study the functoriality of this
new integration  procedure, for example how  to integrate a morphism $A\to B$
of Lie algebroids to a global morphism from a universal stacky
groupoid of $A$ to any stacky groupoid of $B$ (Lie II). The results of
the paper are positive.  We show
\begin{thm}\label{thm:pi2}
The W-groupoid $\cG(A)$ constructed in \cite{tz:arxiv} is
source-2-connected\footnote{Source-2-connectedness means that the source fibres
  have trivial homotopy groups $\pi_{\leq 2}$.}.
\end{thm}

Thus we call $\cG(A)$ the universal W-groupoid of $A$.  Using stacky groupoids, one has one more degree of connectedness. A simple example is
the algebroid 
$A=TS^2$. In this case, $A$ is integrable to the pair groupoid $S^2\times
S^2\rra S^2$, but $\cG(A)= \tilde{S}^2 \times \tilde{S}^2/ B\Z$. Here
$\tilde{S}^2$ is the $B\Z$ gerbe on $S^2$ presented by the action
groupoid $S^3 \times \R \rra S^3$ with $\R$ acting via the
projection $\R \to S^1$ and the usual Hopf $S^1$ action on $S^3$.
Analogously to simply-connected coverings, $\tilde{S}^2$ is the
$\pi_2$-trivial covering of $S^2$ (See Remark \ref{rk:2-conn-cover}). Hence
even with objects as simple as $S^2$ we expect
interesting examples of $\cG(A)$ to appear. Moreover the property
that $\pi_2=0$ might also appeal to symplectic geometers.

But why more connectedness? This is probably best understood via
the following viewpoint elaborated in Section \ref{sec:tgpd}: for
every Lie algebroid $A$, (notice that tangent bundles are Lie
algebroids),  we  associate to $A$
a simplicial set $S(A)=[...S_2(A)\Rrightarrow S_1(A) \Rightarrow
S_0(A)]$ with,
\begin{equation} \label{eq:simp-set}
S_i(A) := \{ \text{Lie algebroid morphisms}\;T\Delta^i \overset{f}{
\to} A, \text{such that $f|_{\text{vertices of}\;\Delta^i}=0$.  }\}.
\end{equation}
Here $\Delta^i$ is the $i$-dimensional standard simplex viewed as
a smooth Riemannian manifold with boundary, hence it is  isomorphic to the
$i$-dimensional closed ball. Then the face and degeneracy
maps are induced by pullbacks of the tangent maps of natural maps $d_k:
\Delta^{i-1} \to \Delta^i$ and $s_k: \Delta^{i} \to \Delta^{i-1}$.
With this language, we also
understand \cite{cafe} and \cite{cf} in a new light:  $S_1(A)$ and
$S_2(A)$ give the space of fields and Hamiltonian
symmetries respectively in the Poisson sigma model \cite{cafe} in the case of
Poisson manifolds; or the space of $A$-paths $P_a A$ and of
$A$-homotopies \cite{cf} respectively in general (see Section \ref{sec:stacky-v}).

The simplicial set $S(A)$ is not entirely unknown to us: in the case of
a Lie
algebra $\g$, let $\Omega^1(\Delta^n, \g)$ be the space of $\g$-valued 1-forms
on $\Delta^n$, then  we have\footnote{More precisely it is $-\alpha$
  the connection 1-form.}, 
\[
\begin{split}
 \hom_{algd}(T\Delta^n, \g)= &\{ \alpha \in \Omega^1(\Delta^n, \g)|
d\alpha = \frac{1}{2}[\alpha , \alpha] \} \\ 
 =& \{ \text{flat connections on the trivial $G$-bundle} \; G\times
 \Delta^n \to \Delta^n\},
\end{split}
\]
where $G$ is a Lie group of
$\g$. 

In fact for  a differential algebra $D$,  
Sullivan \cite{sullivan} constructed a spatial realization of $D$, which is defined as the 
space of differential graded maps
$\hom_{d.g.a.}(D, \Omega^*(\Delta^n))$, and the simplicial set $S(A)$
is a parallel construction on the geometric side. Sullivan's
construction also appears in the work of \v{S}evera \cite{s:funny} to
integrate a non-negatively graded super-manifold with a degree 1 vector
field of square 0 (NQ-manifold), which is the original motivation for our simplicial set
$S(A)$. In fact a Lie algebroid can be viewed as a ``degree 1''
NQ-manifold, and a  Courant algebroid as a
``degree 2'' NQ-manifold \cite{lwx, royt}. Courant algebroids are now widely used in Hitchin's and
Gualtieri's generalized complex geometry \cite{hitchingeneralized, gualtieri}. 
This construction also appears  in the works of  Getzler \cite{getzler} and
Henriques \cite{henriques} to integrate an
$L_\infty$-algebra $L$, where the simplicial set is  $\hom_{d.g.a.} (C^*(L),
\Omega^*(\Delta^n)) $ with $C^*(L) $ the Chevalley-Eilenberg cochains
on $L$.  This simplicial set is further proved to be a Kan simplicial
manifold for all nilpotent $L_\infty$-algebras in \cite{getzler} and for
all $L_\infty$-algebras in \cite{henriques}.

In \cite{z:tgpd-2} we found a 1-1 correspondence between stacky Lie
groupoids and Lie 2-groupoids up to a certain Morita equivalence.
Under this correspondence, $\cG(A)$ corresponds to the {\em 2-truncation} of $S(A)$,
while the usual Lie groupoid of $A$ (when $A$ is integrable)
corresponds to only the {\em 1-truncation} of $S(A)$. It is this  higher
truncation that is the hidden reason for the higher connectedness. Notice that  in the
case of a (finite dimensional) Lie
algebra $\g$, $\cG(\g)$ is the simply connected Lie
group, and this Lie group is always
2-connected.

However it is not obvious that $S(A)$ is a simplicial
manifold let alone a Kan simplicial manifold. Hence its
2-truncation being a Lie 2-groupoid is not immediate but proved
in,

\begin{thm} \label{thm:2-a}
Given a Lie algebroid $A$, the 2-truncation of the
simplicial set $S(A)$,
\[ S_2(A)/ S_3(A) \Rrightarrow S_1(A) \rra
S_0(A), \] is a Lie 2-groupoid that corresponds to the W-groupoid
$\cG(A)$ constructed in \cite{tz} under the correspondence of Theorem 1.4
of \cite{z:tgpd-2}.
\end{thm}
In this theorem $S_{\geq 1}(A)$'s are infinite dimensional spaces. As
a reply to a question of Getzler and Roytenburg, it turns out that
it is not necessary to take everything in the inifite dimensional
space. This Lie 2-groupoid is Morita equivalent to a finite
dimensional Lie 2-groupoid, $E\Rrightarrow P \rra M$, arising in the form of local Lie
groupoids \`a la Pradines,  where $\dim P= \dim
A$ and $\dim E=2\dim A -\dim M$ (See \eqref{eq:2gpd-ft}). 

Finally, after finding  this universal W-groupoid $\cG(A)$,
we have the expected
\begin{thm} [Lie II for Lie algebroids] \label{thm:lieii}
Let $\phi$ be a morphism of Lie algebroids $A \to B$, $\cG$ a
W-groupoid whose algebroid is $B$. Then up to 2-morphisms, there
exists a unique morphism $\Phi$ of W-groupoids $\cG(A)\to \cG$ such
that $\Phi$ induces the Lie algebroid morphism $\phi: A\to B$.
\end{thm}

\noindent \noindent {\bf Acknowledgments:} I thank Alan
Weinstein for his continuing guidance and support.  I thank Andre
Henriques for much helpful discussion and for communicating
preliminary versions of his
work on string groups. I am also very grateful to 
Ezra Getzler, Kirill Mackenzie, Pavel \v{S}evera for very helpful
discussions. Lastly, I thank the referee for very useful suggestions to improve the paper.

\section{Technical tools} 
\label{sec:tt} In this section we recall some basic knowledge on the technical tools that we need later in the paper.

\subsection{Groupoids, stacks and higher groupoids}
We first recall necessary knowledge in the theory of Lie groupoids and
differentiable stacks \cite{bx,
  metzler, pronk}  \cite{m-orbi} (in the context of orbifolds). {\em
  An H.S. (Hilsum-Skandalis) morphism} between two Lie groupoids
$K:= K_1 \rra K_0$ and $K':= K'_1 \rra K'_0$ is given by a
right-principal bibundle $E$, which is called the {\em H.S.
bibundle},
\[
\xymatrix{ 
K_1 \ar[d] \ar@<-1ex>[d] & E \ar[dl]_{J_l} \ar[dr]^{J_r} & K'_1 \ar[d]
\ar@<-1ex>[d] \\
K_0 & & K'_0.
}
\]
More precisely, a {\em bibundle} $E$ is a manifold on which both $K$ and
$K'$ act compatibly
\[\Phi_l : K_1 \times_{\bs, K_0, J_l} E \to E, \quad \Phi_r:
E\times_{J_r, K'_0, \bt} K'_1 \to E, \]
via the  moment maps $J_l$ and $J_r$.  $E$ is furthermore
said to be {\em right-principal}, if the right action of $K'$ is principal, that is,
$J_l$ is a surjective submersion and $id \times \Phi_r: E\times_{J_r,
  K'_0, \bt}K'_1 \cong
E\times_{J_l, K_0, J_l} E $ is an isomorphism. If $E$ is also
left-principal (defined in a similar way), $E$ is called a {\em Morita
  bibundle} and it gives a {\em Morita equivalence} between $K$ and
$K'$. Then there is a dictionary between differentiable stacks and Lie
groupoids: 
\begin{center}
\begin{tabular} {ccc}
differentiable stacks & $\leftrightarrow$ & Lie groupoids \\
morphisms &  $\leftrightarrow$  & H.S. morphisms \\
isomorphisms &  $\leftrightarrow$ & Morita equivalences.
\end{tabular}
\end{center} 
Given a differentiable stack $\cX$, a Lie groupoid $X:=X_1 \rra X_0$ corresponding to it
is called a {\em groupoid  presentation} of $\cX$, and $X_0$ is called a {\em
  chart} of $\cX$. 

There are also various version of higher Lie groupoids. One uniform
way to describe Lie $n$-groupoids is via its nerve, by requiring the nerve to
be a simplicial manifold whose homotopy groups are trivial above
$\pi_n$. We refer the reader to the introduction of \cite{z:tgpd-2} and
the references therein for a general discussion of these matters, and only
recall briefly the definitions we need here.

A simplicial set (resp. manifold) $X$ is made up by sets (resp.
manifolds) $X_n$ and structure maps
\begin{equation}\label{eq:fd}
 d^n_i: X_n \to X_{n-1} \;\text{(face maps)}\quad s^n_i: X_n \to X_{n+1} \; \text{(degeneracy maps)},\;\; i\in \{0, 1, 2,..., n\} \end{equation}
that satisfy suitable coherence conditions. The first two examples
are the simplicial $m$-simplex $\Delta[m]$ and the horn
$\Lambda[m,j]$, defined as
\begin{equation}\label{eq:simplex-horn}
\begin{split}
(\Delta[m])_n & = \{ f: (0,1,...,n) \to (0,1,..., m)| f(i)\leq
f(j),
\forall i \leq j\}, \\
(\Lambda[m,j])_n & = \{ f\in (\Delta[m])_n| \{0,...,j-1,j+1,...,m\}
\nsubseteq \{ f(0),..., f(n)\} \}.
\end{split}
\end{equation}

\begin{ep}\label{ep:nerve}
Given a Lie groupoid $G:= G_1 \rra G_0 $, we complete it to a simplicial
manifold $NG$, with 
\[ NG_n = G_1 \times_{\bs, G_0, \bt} G_1 \times_{ \bs, G_0, \bt} \dots
\times_{ \bs, G_0, \bt} G_1, \quad \text{$n$ copies of $G_1$}, \] for
$n\ge 1$ and $NG_0=G_0$. The face and degeneracy maps are
\begin{eqnarray*}
& d_k(g_1,{\dots},g_n ) =   & s_k(g_1, \dots, g_n) =\\
& \begin{cases}
  (g_2, {\dots}, g_n) & k=0 \\
  (g_1, {\dots}, g_k g_{k+1}, {\dots}, g_n) & 0<k<n, \\
  (g_1, {\dots}, g_{n-1})  & k=n,
  \end{cases}  & \begin{cases} (1_{\bt(g_1)}, g_1,
    \dots, g_n) & k=0 \\
(g_1, \dots, g_k,  1_{\bt(g_{k+1})}, g_{k+1}, \dots, g_n) & 0<k<n, \\
(g_1, \dots, g_n, 1_{\bs(g_n)}) & k=n.
\end{cases} 
\end{eqnarray*}This construction is known as the {\em nerve of a Lie
groupoid}.
\end{ep}
A simplicial set $X$ is {\em Kan} if any map from the horn
$\Lambda[m,j]$ to $X$ ($m\ge 1$, $j=1,..,m$), extends to a map from
$\Delta[m]$ to $X$. In the language of groupoids, the Kan condition
corresponds to the possibility of composing various morphisms. In
an $n$-groupoid, the only well defined composition law is the one
for $n$-morphisms. This motivates the following definition.

\begin{defi}\cite{henriques} \label{def:defngroupoid} A {\em Lie} $n${\em-groupoid} $X$
($n\in\N \cup \infty$) is a simplicial manifold that satisfies the conditions
$Kan(m,j)$ $\forall m\ge 1$, $0\le j\le m$ and $Kan!(m,j)$
$\forall m> n$, $0\le j\le m$, where: \vspace{.3cm}

\noindent
\begin{tabular}{p{1.6cm}p{12.6cm}}
$Kan(m,j)$:& The restriction map
$\hom(\Delta[m],X)\to\hom(\Lambda[m,j],X)$ is a
surjective submersion.\\
$Kan!(m,j)$:& The restriction map
$\hom(\Delta[m],X)\to\hom(\Lambda[m,j],X)$ is a
diffeomorphism.\\
\end{tabular} \vspace{.3cm} \\
When $n=\infty$, a Lie $\infty$-groupoid is also called a {\em Kan
  simplicial manifold}. 
\end{defi}

\begin{remark}\label{rk:finite} A Lie $n$-groupoid $X$ is determined by
its first $(n+1)$-layers $X_0$, $X_1,\dots, X_n$ and some structure
  maps. For example a Lie 1-groupoid is exactly determined by a Lie
  groupoid structure on $X_1\rra X_0$. In fact a Lie 1-groupoid is the
  nerve of a Lie groupoid. A Lie 2-groupoid is exactly
  determined by $X_2 \Rrightarrow X_2 \rra X_0$ together with a sort of
  3-multiplications and face and degeneracy maps satisfying certain
compatibility conditions. This is made precise in \cite[Section 2]{z:tgpd-2}.
Hence in this paper, we often write only the first three layers of a
Lie 2-groupoid. 
\end{remark}

\begin{remark}\label{rk:sheaf}
It is not clear that $\hom(\Lambda[m,j], X)$ is a manifold. In
general, as in \cite[Section 2]{henriques}, we sometimes talk about the
limit (for example fibre-product) of a diagram  in $\cC$, the category of
differential manifolds, before knowing its existence. For this purpose, we use the Yoneda functor
\[ 
\begin{split} 
yon: \cC & \to \{ \text{Sheaves on $\cC$} \} \\ 
 X & \mapsto (T\mapsto \hom(T, X))
\end{split}
\]
to embed $\cC$ to the category of sheaves on $\cC$. Hence a limit of
objects of $\cC$ can always be viewed as the limit of the
corresponding representable sheaves using $yon$. The limit sheaf is
representable if and only if the original diagram has a limit in
$\cC$. Viewed as a sheaf, $\hom(\Lambda[m,j], X)$ is
representable by \cite[Lemma
2.4]{henriques}, thus it is a manifold always. 
\end{remark}

For our convenience, (since in the end we will apply Kan replacement to local
Lie groupoids), we need simplicial manifolds with certain properties, 
described  below,
\begin{enumerate}
\item \label{itm:a} {\bf Property A}: for all $0\le j \le m$, the restriction map
$\hom(\Delta[m],X)\to\hom(\Lambda[m,j],X)$ is a submersion;
\item {\bf Property B}: there are isomorphisms 
\[ \hom(\Lambda[2, 0], X)\cong \hom(\Lambda[2,1], X) \cong \hom
(\Lambda[2,2], X) ,\]  
which are compatible with the face maps $
\hom(\Delta[2], X) \to \hom(\Lambda[2, j], X)$. 
\end{enumerate}
 
\begin{ep}\label{ep:localgpd}
A {\em local Lie groupoid} $G^{loc}:=(U, V) \rra M$ resembles a Lie
groupoid, with the difference that the multiplication is only defined on $V \times_{M} V$
where $V$ is an open tubular neighborhood of the identity section
$M\to U$. We also require here that the inverse map $i$ maps $V$ into $V$.  This  gives us a simplicial manifold
$NG^{loc}$ (the nerve of local Lie groupoid), with $NG^{loc}_0=M$, $ NG^{loc}_1 = V$,
\[ NG^{loc}_n= \{ (g_1, \dots, g_n)\in V\times_M V \times_M \dots \times_M V| g_i\cdot \dots \cdot g_{j} \in V, \forall 1\le i \le j \le n \}, \]
and face and degeneracy maps similar to the ones in the nerve of a groupoid (Example \ref{ep:nerve}).
Then we have 
\begin{lemma}
Constructed as above, $NG^{loc}$ is a simplicial manifold satisfying properties A and B. 
\end{lemma}
\begin{proof}
First of all, we notice that the multiplication $m: V \times_M V \to
U$ is a submersion: it is easy to argue that the groupoid
multiplication is surjective, by using inverse elements. By lifting
this argument to the tangent groupoid, we see that the multiplication
map is a submersion. Since being a submersion is a local property, this argument applies to local Lie groupoids. Thus $NG^{loc}_2=m^{-1}(V)$ is a manifold. The natural embedding $NG^{loc}_2 \to \hom(\Lambda[2,1], NG^{loc}) = V\times_M V$ is a submersion. Moreover, with the inverse map $i: V\to V$, it is easy to see that $NG^{loc}$ satisfies Property B. Thus $NG^{loc}$ satisfies Property A for $m=1,2$.

It is not hard to see that $NG^{loc}_n = \hom(\Lambda[n, j], NG^{loc})$  for all $0\le j \le n$; for example,
\[
\begin{split}
&\hom(\Lambda[3,1], NG^{loc})=\{(g_1, g_2, g_3)|(g_1, g_2) \in NG^{loc}_2, (g_1, g_2 g_3) \in NG^{loc}_2, (g_2, g_3) \in NG^{loc}_2\} \\
= &\{(g_1, g_2, g_3)|g_i\cdot \dots \cdot g_{j} \in V, \forall 1\le i \le j \le 3\} \\
= & NG^{loc}_3.
\end{split}
\] 

The induction in  \cite[Lemma 2.4]{henriques} shows that  if we know
that $NG^{loc}$ satisfies Property A for $m\le n-1$ (thus naturally
$NG^{loc}_{m}$ and $\hom(\Lambda[m, j], NG^{loc})$ are representable
for $m\le n-1$), then $\hom(\Lambda[n,j], NG^{loc})$ is representable. Thus $NG^{loc}_n = \hom(\Lambda[n, j], NG^{loc})$ is a manifold. 

Thus  $NG^{loc}$ is actually a Lie 2-groupoid (the point is that it is a simplicial manifold) satisfying properties A and B. 
\end{proof}
\end{ep}

There is another sort of higher groupoid, whose space of arrows is a differentiable stack
$\cG$ and whose space of object is a manifold $M$.  $\cG \rra M$ satisfies the usual
groupoid axioms up to some 2-morphisms, which in turn satisfy some
higher coherences.  We call such an object a {\em stacky Lie (SLie for short)
groupoid}. When the stack is
further \'etale, we call it a {\em W-groupoid} as in the  introduction. We refer the reader
to \cite{z:tgpd-2} for a complete definition and the following result: There is a
1-1 correspondence of Lie 2-groupoids and SLie groupoids up to a
certain equivalence. Given an SLie groupoid $\cG \rra M$
with a groupoid presentation $G:=(G_1\rra G_0)$ of $\cG$, its associated Lie
2-groupoid $Y$ is given by
\begin{equation} \label{eq:s-2}
Y_0 =M , \quad Y_1=G_0, \quad Y_2=E_m, 
\end{equation}
where $E_m$ is the bibundle representing the multiplication $m:
\cG\times_{M} \cG \to \cG$ of the SLie groupoid. Conversely,
given a Lie 2-groupoid $Y$, let  $b(Y_2):= d_2^{-1}(s_0(Y_0))$ be the
set of bigons in $Y_2$. Then $d_0, d_1: b(Y_2) \rra Y_1$ is a Lie
groupoid, and  the stack $\cG$ presented by $b(Y_2) \rra Y_1$ is an
SLie groupoid over $Y_0$. 

\subsection{ Hypercovers, generalized morphisms, Morita equivalences,
  and truncations}

\begin{defi}\label{defequivalence}
A strict map $f:Z\to X$ of Lie $n$-groupoids is a {\em hypercover} if the natural
map from $Z_k=\hom(\D^k,Z)$ to the pull-back
\begin{equation}\label{eq:hc-rhs}
\pb(\pD[k],Z, \D[k], X) =PB(\hom(
\partial \Delta[k] , Z)\ra \hom(\partial \Delta[k], X) \la X_k)
\end{equation}
is a surjective submersion for $0\le k\le n-1$ and an
isomorphism\footnote{When $n=\infty$, namely in the case of Kan
  simplicial manifolds, the requirement of isomorphism is empty.} for $k=n$. 
\end{defi}

\begin{defi}\label{defi:gen-morp} A {\em generalized morphism} of Lie $n$-groupoids from
  $X$ to $Y$ consists of a {\em zig-zag}  of strict maps
$X\stackrel{\sim}{\leftarrow}Z\to Y$, where the map $Z\stackrel{\sim}{\to} X$ is
an hypercover. 
\end{defi}

\begin{defi}
A {\em $2$-morphism}\footnote{This is not the most general notion of
  2-morphism, however it is sufficient for the purposes of this paper.}
between two generalized morphisms of Lie $n$-groupoids, 
$X\stackrel{\sim}{\leftarrow}Z\to Y$ and
$X\stackrel{\sim}{\leftarrow}Z'\to Y$, is an isomorphism of Lie
$n$-groupoids $Z\cong Z'$ which commutes with the maps $Z, Z' \to X$ and
$Z, Z' \to Y $.
\end{defi}

\begin{defi}\label{defi:m-equi-2gpd} Two Lie $n$-groupoids $X$ and
  $Y$ are
{\em Morita equivalent} if there is a Lie $n$-groupoid $Z$ such that
both of the maps $ X\stackrel{\sim}{\leftarrow} Z
\stackrel{\sim}{\rightarrow} Y$ are hypercovers.
\end{defi}

We denote a hypercover as $X \stackrel{\sim}{\to} Y$,  a
generalized morphism as $\xymatrix@C=.45cm{X \ar@{.>}[r] & Y}$, and a Morita
equivalence as $X \stackrel{\sim}{\leftrightarrow}  Y$. 

Hypercovers of Lie $n$-groupoids may also be understood as higher
analogues of pull-backs of Lie groupoids via a surjective submersion. Let $X$ be a 2-groupoid, 
let $Z_1\rra Z_0$ be two
manifolds with structure maps as in \eqref{eq:fd} up to
the level $n = 1$, 
and let $f_n: Z_n \to X_n$ be a map
preserving the
structure maps $d^{n}_k$'s and $s^{n-1}_k$'s for $n\leq 1$. Then
$\hom(
\partial \Delta[n] , Z)$ 
is defined
for $n\leq 1$. We further
suppose that
$f_0: Z_0 \twoheadrightarrow X_0$ is a surjective submersion
(hence $Z_0\times Z_0
\times_{X_0\times X_0} X_1$ is a manifold) and do the same for $Z_1\thra Z_0\times
Z_0 \times_{X_0\times X_0} X_1$.
That is to say, we suppose that the induced
maps from $Z_k$ to the pull-backs $\hom(
\partial \Delta[k] , Z)\times_{ \hom(\partial \Delta[k], X)}  X_k$
are surjective submersions for $k= 0, 1$. Then we form
\[Z_2= \hom(
\partial \Delta[2] , Z)\times_{ \hom(\partial \Delta[2], X)}  X_2,\]
which is a manifold (see \cite[Lemma 2.4]{z:tgpd-2}). 

Moreover there are maps $d^2_i: Z_2\to Z_1$ induced by the natural
projections $\hom(\partial \Delta[2] , Z)\to Z_1$, and maps $s^1_i: Z_1 \to
Z_2$ 
defined by
\[ s^1_0(h)=(h,h,s^0_0(d^1_0(h)),s^1_0(f_1(h))), \quad
s^1_1(h)=(s^0_0(d^1_1(h)),h,h,s^1_1(f_1(h))).
\]
Similarly we define higher $Z_i$'s inductively. As shown in
\cite[Section 2.2]{z:tgpd-2},  $Z_2
\Rrightarrow Z_1 \rra Z_0$ is a Lie $2$-groupoid and we call it the
{\em pull-back Lie 2-groupoid}
of $X$ by $f$. 
Then $f: Z\to X$ is a hypercover with $f_2$ the natural projection
$f_2: Z_2 \to X_2$. In particular, $Z$ is Morita equivalent to
$X$. Sometimes $f$ does not make the induced
maps from $Z_k$ to $\hom(
\partial \Delta[k] , Z)\times_{ \hom(\partial \Delta[k], X)}  X_k$
surjective submersions for $k= 0, 1$, but $Z_2$ is still a
manifold. In this case, we can still construct a pull-back
Lie $2$-groupoid as above. However, the pullback is not obviously
Morita equivalent to $X$ any more. 
Nevertheless, the following lemma, that we will use later, says that
in a certain case this still holds. 

\begin{lemma}\label{lemma:me}
Given a simplicial morphism $f:K \to Y$ of two Lie 2-groupoids,
if 
\begin{enumerate}
\item\label{itm:1} $K_0\cong Y_0$, 
\item\label{itm:2} there are a manifold $Z_1$, and  surjective submersions $K_1 \xleftarrow{g_1} Z_1 \xrightarrow{h_1} Y_1$, 
\item\label{itm:additional}    there exists  a morphism $\alpha:  Z_1 \to b(Y_2)$ such that $f_1 \circ
  g_1= \alpha \cdot h_1$,   
\item \label{itm:iso} the natural map $K_2 \to \hom(\partial \Delta[2]\to \Delta[2], K
  \to Y)$ is an isomorphism,
\end{enumerate}
then $K$ is Morita equivalent to $Y$.
\end{lemma}
\begin{proof}
Take $Z_0=K_0=Y_0$,   $g_0=h_0=id$. 
We pullback $K$ and $Y$ via $g_1$
and $h_1$ to $Z_1$. To show $K$ and $Y$ are Morita equivalent, we only
need to show that the map induced by $f$
\begin{equation}\label{eq:induced-by-f}
 \hom(\partial \Delta[2]\to \Delta[2], Z
  \to Y) \to \hom(\partial \Delta[2]\to \Delta[2], Z
  \to K), 
\end{equation}
is an isomorphism. 
Then \eqref{eq:induced-by-f} follows from items \ref{itm:additional} and \ref{itm:iso} :
\[
\begin{split}
\text{RHS}= & \hom(\partial\Delta[2], Z)
  \times_{h_1, \hom(\partial\Delta[2], Y)} \hom(\Delta[2], Y) \\
=&\hom(\partial\Delta[2], Z)
  \times_{f_1 \circ g_1, \hom(\partial\Delta[2], Y)} \hom(\Delta[2], Y) \\
=&\hom(\partial\Delta[2], Z)
  \times_{g_1, \hom(\partial\Delta[2], K)} \hom(\partial\Delta[2], K)
  \times_{f_1, \hom(\partial\Delta[2], Y)} \hom(\Delta[2], Y)\\
=&\hom(\partial\Delta[2], Z)
  \times_{\hom(\partial\Delta[2], K)} \hom(\Delta[2], K)=\text{LHS}.
\end{split}
\]
\end{proof}

\begin{remark}\label{rk:needed} Given a simplicial morphism $f:K \to Y$ of  Lie 2-groupoids, if 
\begin{enumerate}
\item\label{itm:1'} $K_0\cong Y_0$  and there is a submanifold $X_1$ of
  $K_1$ such that $X_1\xrightarrow{f_1|_{X_1}}
  Y_1$ is a submersion, 
\item\label{itm:2'} the composed map $h'_1: \hom(\Lambda[2,1], K)\times_{\hom(\Lambda[2,1],Y)}
Y_2 \to Y_2\xrightarrow{d_1} Y_1$ is a surjective submersion,
\item\label{itm:additional'}    there exist a morphism $\mu: \hom(\Lambda[2,1],
 K) \to K_1$ and a morphism $\alpha':  \hom(\Lambda[2,1],
 K)\times_{\hom(\Lambda[2,1],Y)} 
Y_2 \to b(Y_2)$ such that the composed map $g'_1:   \hom(\Lambda[2,1],
 K)\times_{\hom(\Lambda[2,1],Y)} 
Y_2  \to \hom(\Lambda[2,1],
  K) \xrightarrow{\mu} K_1$ makes $f_1 \circ
  g'_1= \alpha \cdot h'_1$,  
\item \label{itm:surjective'} the natural map  $g_1:=g'_1\sqcup pr_{X_1}: Z_1:= \hom(\Lambda[2,1], K) \times_{\hom(\Lambda[2,1],Y)} 
Y_2 \sqcup X_1 \times_{Y_1, d_0} b(Y_2) \to Y_1$ is a surjective submersion, 
\item \label{itm:iso'} the natural map $K_2 \to \hom(\partial \Delta[2]\to \Delta[2], K
  \to Y)$ is an isomorphism,
\end{enumerate}
then $Z_1$, $g_1$, $h_1:=h'_1 \sqcup
  d_1 \circ pr_{b(Y_2)}$ and an extension of $\alpha'$ to $\alpha: Z_1
  \to b(Y_2)$ obtained  by addition of
$\alpha=pr_{b(Y_2)}: X_1 \times_{Y_1, d_0} b(Y_2) \to b(Y_2)$ satisfy the conditions required by the above lemma.
\end{remark}

\begin{defi} \label{defi:truncation}
Given a Kan simplicial manifold $X$, we can perform a certain
operation $\tau_n$ called $n${\em -truncation} (it is denoted by
$\tau_{\le n} $ in \cite[Section 3]{henriques}) for $n\in \Z^{\ge 1}$, 
\[  \tau_n(X)_k = X_k, \forall k\le n-1, \quad \tau_n(X)_k =
X_k/\sim_k, \forall k \ge n, \]
where two elements $x\sim_k y $ in $X_k$ if they are simplicially homotopic relative to their $(n-1)$-skeleton\footnote{When $k=n$, this means explicitly that $d_i x = d_i y$, $0\le i \le k$,
  and there exists  $z\in X_{k+1}$ such that $d_k(z)=x, d_{k+1}(z)=y$,
and $d_i z = s_{k-1} d_i x = s_{k-1} d_i y$, $0\le i < k$. }. 
\end{defi}


\subsection{Kan replacement}

Here we recall the construction of a canonical functor $Kan$ from the category of
simplicial manifolds with  properties A and B (they are called invertible local Kan simplicial manifolds in \cite{z:kan}) to the category of weak Kan
simplicial manifolds. 

Let $J$ be the set of inclusions as below, 
\begin{equation}\label{eq:j} J:= \{ \Lambda[k, j] \to \Delta[k]: 0\le j \le k\ge 3,  \} \cup \{ \Lambda[2, 1]
\to \Delta[2] \}. \end{equation}
Given a simplicial manifold with property A, we construct a series of simplicial manifolds
\begin{equation}\label{eq:series}
X=X^0 \to X^1 \to X^2 \to \dots \to X^\beta \to \dots  
\end{equation}
by an inductive pushout:
\begin{equation}\label{diag:x}
\begin{diagram}
\sqcup_{(\Lambda[k, j] \to \Delta[k]) \in J} \Lambda[k, j] \times
\hom(\Lambda[k, j], X^\beta)   &    \rTo    &  X^\beta \\
 \dTo       &   & \dTo \\
\sqcup_{(\Lambda[k, j] \to \Delta[k]) \in J} \Delta[k]
\times\hom(\Lambda[k, j], X^\beta)  &    \rTo    &  {\NWpbk} X^{\beta+1}&.
\end{diagram}
\end{equation}
Then we let $Kan(X) = \colim_{\beta \in \N} X^\beta$. The construction
of $Kan(X)$ makes it clear that the morphism  $ X \to Kan(X) $  is
functorial. However $Kan(X)$ is only a weak Kan simplicial manifold,
whose definition is omitted here since we will not use it in this paper. But we still have
\begin{thm}\label{lemma:kankan} \cite[Thm. 3.6]{z:kan} If $W$ is Lie 2-groupoid, then $\tau_2(Kan(W))$ is a Lie 2-groupoid which is Morita equivalent to $W$. 
\end{thm}

Notice that given any finite simplicial set $A$ (both
$\Lambda[k,j]$ and $\Delta[k]$ are such), the natural map of sets is an
isomorphism
\begin{equation}\label{eq:iso-colim}\colim_{\beta} \hom(A, X^\beta)
  \xrightarrow{\simeq} \hom(A, Kan(X)). 
\end{equation}

We make some calculation for Kan replacement:

First of all, $X_0 = X^1_0 = X^2_0 = \dots
= Kan(X)_0$, and 
\begin{equation}\label{eq:x1}
\begin{split}
X^1_1 &=  X_1 \sqcup (X_1 \times_{X_0} X_1) \\
X^2_1 &= X^1_1 \sqcup X^1_1 \times_{X_0} X^1_1 \\
      &= X_1^1 \sqcup \big( X_1 \times_{X_0} X_1
\sqcup X_1 \times_{X_0} (X_1 \times_{X_0} X_1) \\
&\sqcup (X_1 \times_{X_0} X_1)  \times_{X_0} X_1 \sqcup  (X_1
\times_{X_0} X_1)\times_{X_0}(X_1 \times_{X_0} X_1) \big) \\
& \vdots \\
 Kan(X)_1 &= X_1 \sqcup ( X_1 \times_{X_0} X_1) \sqcup (X^1_1
 \times_{X_0} X^1_1 ) \sqcup  (X^2_1 \times_{X_0} X^2_1)  \dots,  
\end{split}
\end{equation}
which can be represented by the following picture: \\
\psfrag{...}{$\dots$} \psfrag{X1}{$Kan(X)_1:$}
\centerline{\epsfig{file=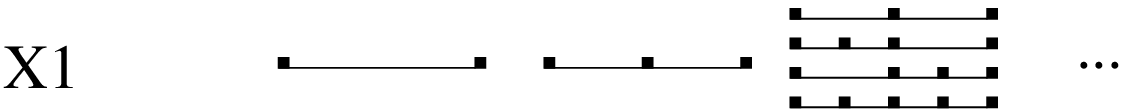, height=1cm, width=8cm}}

A calculation shows that
\[
\begin{split}
X^1_2 &= X_2 \sqcup X_1 \times_{X_0} X_1 \sqcup X_1 \times_{X_0} X_1
\sqcup X_1 \times_{X_0} X_1 \\
& \sqcup (\sqcup_{j=0}^3\hom(\L[3, j], X) \\
X^2_2 &= X^1_2 \sqcup X^1_1 \times_{X_0} X^1_1 \sqcup X^1_1 \times_{X_0} X^1_1
\sqcup X^1_1 \times_{X_0} X^1_1 \\
& \sqcup (\sqcup_{j=0}^3\hom(\L[3, j], X^1) \\ 
&\vdots
\end{split}
\]
Inside $X^1_2$, there are three copies of $X_1\times_{X_0} X_1$. The
first is an artificial filling of the horn $X_1\times_{X_0} X_1$, and
the second two are images of degeneracies of $X_1\times_{X_0} X_1$ in
$X^1_1$. The same for $X^2_2$, etc.  We represent an element in $X^1_2$
as\\
\psfrag{X2}{$X_2:$}\psfrag{X11}{$X_1
  \times_{X_0}X_1:$}\psfrag{hom}{$\hom(\L(3, j), X):$}
\psfrag{...4 of them}{$\dots$ 4 such} \psfrag{,}{,}
\begin{equation}\label{pic:x12}
\epsfig{file=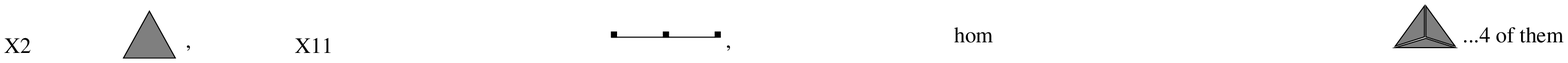, height=.5cm, width=12cm }
\end{equation}
plus those degenerate ones in the other two copies of  $X_1
\times_{X_0} X_1$. Thus further we represent a non-degenerate element
in $X^2_2$ as 

\vspace{.5cm}
\psfrag{X21}{$X_2^1$: described as above}
\psfrag{X11}{$X_1^1\times_{X_0}X_1^1$:} \psfrag{X3j}{$\hom(\L[3, j],
  X^1)$:} \psfrag{... 7 of them}{$\dots$}
\begin{equation}\label{pic:x22}
\epsfig{file=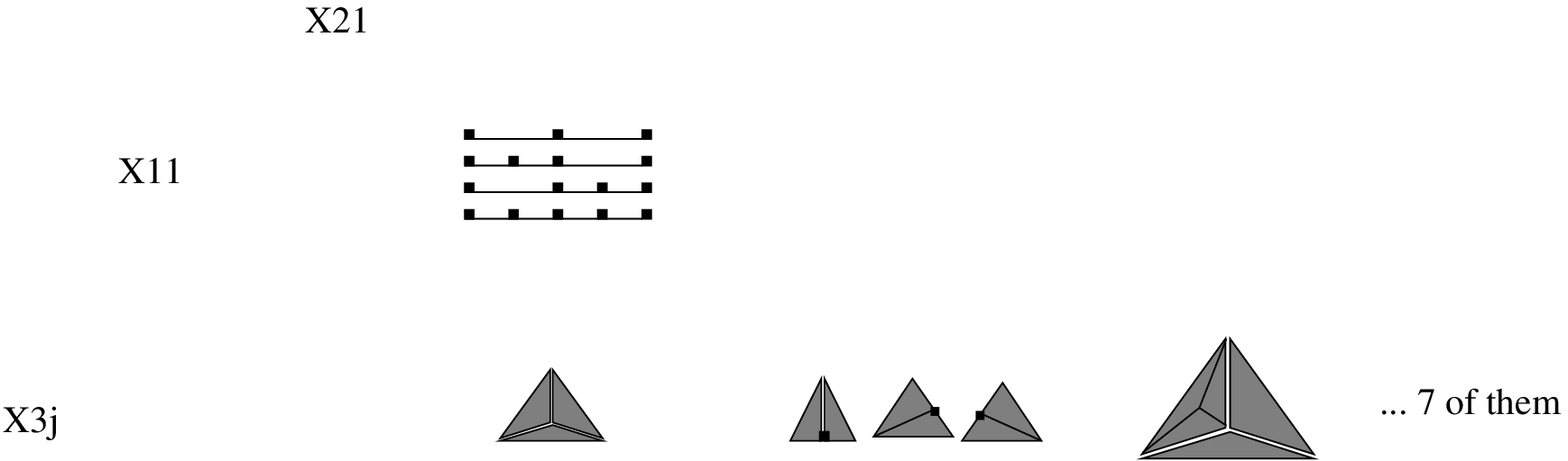, height=3cm, width=12cm}
\end{equation}

\section{The universal groupoid of an algebroid $A$} \label{sec:tgpd} \label{sec:wgpd} 

The first construction of universal groupoid of a Lie algebroid $A$
(even though not the one we use here) is contained in \cite{cf}, where
a topological groupoid (called Weinstein groupoid therein) is
constructed via a certain path method which we will recall soon. This
topological groupoid can be enriched with a differentiable stack
structure \cite{tz} and the latter is  the universal groupoid we use
in this paper.  We are going to show that this path method can be extended to higher levels, known  as Sullivan's construction. This gives a simplicial version of this universal stacky groupoid. Then we will also give a Lie 2-groupoid version and show the relation between all these versions.  

\subsection{The stacky groupoid version}\label{sec:stacky-v}
 We first recall \cite{cf} that an
{\em $A$-path} of a Lie algebroid $A$ is a $C^1$ path  $a(t): \Delta^1
\to A$ with base map
$\gamma(t): \Delta^1 \to M$ satisfying
\begin{equation}\label{eq:apath}
\rho(a(t))=\dot{\gamma}(t),
\end{equation} where $\dot{\square}$ denotes the derivative of $t$.
If it further satisfies boundary condition
$a(0)=a(1)=\dot{a}(0)=\dot{a}(1)=0$, we call it an
$A_0$-path as in \cite{tz}. This notation is  useful for technical reasons (see Theorem \ref{thm:2-a'}).
We call $P_a A$ and $P_0 A$ the space of $A$-paths and $A_0$-paths
respectively. By \cite[Lemma 4.6]{cf} and a slightly modified argument (see also \cite[Section 2]{tz:arxiv}),  $P_a A$ is a Banach manifold and  $P_0 A$ is a Banach submanifold of $P_a A$.

\begin{remark}\label{rk:hilbert}
Since we choose $C^1$-paths, our Banach manifolds in this paper are all Hilbert manifolds, namely manifolds based on Hilbert open charts. In fact, as long as we take $C^r$-morphisms, and $r$ is finite, this holds always true \cite[Chap.3.1]{Brylinski}. However we keep on the terminology of Banach manifolds to match with  the references such as \cite{henriques, lang}. 
\end{remark}

There is an equivalence relation in $P_a A$, called $A$-homotopy \cite{cf}.
\begin{defi}
Let $a(t, s)$ be a $C^1$-family of $A$-paths. Assume that the base paths $\gamma(t,
s):=\rho\circ a(t,s)$ have fixed end points. Fix a connection\footnote{Two $A$-paths being $A$-homotopic is in the end independent of the choice of connections.} $\nabla$ on $A$. Two $A$-paths $a_0$ and $a_1$ are $A$-homotopic if there are $C^1$-morphisms $a(t, s)$, $b(t, s)$: $\Delta^1 \times \Delta^1 \to A$ satisfying 
\begin{equation}\label{eq:hom}
\partial_t b-\partial_s a=T_\nabla(a,b),\quad b(0,s)=0, \quad b(1, s)=0
\end{equation} Here $T_\nabla$ is the torsion of the connection defined
by $T_\nabla(\alpha, \beta)=\nabla_{\rho(\beta)}\alpha-\nabla_{\rho(\alpha)}\beta +
[\alpha,\beta].$ The restriction of an $A$-homotopy on $P_0 A$ gives an $A$-homotopy of $A_0$-paths.
\end{defi}
\begin{remark}
 This definition is only apparently different from the one in \cite{cf} by Lemma 1.5 (iii) therein. We need this version since it coincides with the simplicial picture. 
\end{remark}

This flow of $A_0$-paths $a(\epsilon, t)$ generates a foliation
$\cF$ of codimension $\dim
A$ on $P_0A$.  The $A_0$-path
space is a Banach submanifold of the $A$-path space and $\cF$ is the
restricted foliation of the foliation defined in 
\cite[Section 4]{cf}. Sometimes to avoid
dealing with infinite dimensional issues, we take an open cover $U_i$ of
$P_0 A$, and take $G_0(A)$  the
disjoint union of slices $P_i\subset U_i$  that are transversal to
the foliation $\cF$.  As the classic foliation theory tells us (see for example \cite[Example 5.9]{moerdijk}),  $G_0(A)$ is a smooth immersion submanifold
of $P_0 A$, and $\cF$
induces an \'etale foliation $\cF_{\acute{e}}$ on $G_0(A)$.  Also recall that the monodromy groupoid assigned to a foliation is the
groupoid whose space of objects is the manifold itself and whose space
of arrows is given by leafwise paths modding out leafwise homotopies.  We take
the monodromy groupoid 
$G(A):=G_1(A)\rra G_0(A)$ of $\cF_{\acute{e}}$.  It is  a finite dimensional \'etale
Lie groupoid
and  it is  the
pull-back groupoid of the monodromy groupoid  $Mon(P_0 A):=Mon(P_0 A)_1 \rra P_0 A$ via $G_0(A) \to
P_0 A$. Moreover $Mon(P_0 A)_1 \times_{P_0 A} G_0(A) \to P_0A$ is surjective and \'etale, thus $G(A)$ and $Mon(P_0 A)$ are Morita equivalent.


There are  two maps
$G_0(A) \rra M$ assigning to a path its two end points. 
The quotient stack
$\cG(A)_1:=[G_0(A)/G_1(A)]$, which is presented by $G(A)$,  is 
an \'etale differentiable stack. The two maps from $G_0(A)$ to $M$
descend to the quotient, giving two maps $\bs, \bt:\cG(A)\to M$. We define a
multiplication $m:\cG(A)_1 \times_{\bs,M,\bt}\cG(A)_1\to \cG(A)_1$  by
concatenation $\odot$ of paths  (here we need to use the
  Banach chart $P_0 A$), namely,
\begin{equation}\label{eq:con-2}
(a\odot b)(t)
=\begin{cases}  a(2t)\quad \text{when}\; t \in [0,1/2], \\
b(2t-1)\; \text{ when}\;
t \in [1/2, 1];
\end{cases}
\end{equation}
We define an inverse
$i:\cG(A)_1\to\cG(A)_1$, by reversing
the orientation of a path; we  define an
identity section $e: M\to \cG(A)_1$ by considering constant paths. These maps are defined in
detail in \cite{tz}. There,  we prove that these structure maps make
$\cG(A):=\cG(A)_1\rra M$ into a W-groupoid.

\subsection{The simplicial version}\label{sec:simp-v}
Now to a Lie algebroid $A$ over a manifold $M$, we  associate the
simplicial set $S(A)$ as in \eqref{eq:simp-set} in the introduction. The first three
layers of $S(A)$ are actually familiar to us:
\begin{itemize}
\item  it is easy to check that $S_0=M$;
\item  $S_1$ is exactly the $A_0$-path space $P_0 A$
since a map $T \Delta^1 \to A$ can be written as $a(t) dt$ with
base map $\gamma(t): \Delta^1 \to M$.  Its being a Lie algebroid
map is equivalent to $\rho(a(t))=\frac{d}{dt}\gamma(t)$ since the
Lie bracket of $T \Delta^1$ is trivial and its anchor is the
identity;
\item bigons in $S_2$ are exactly the $A$-homotopies in $P_0 A$ since a bigon $f: T (d^2_2)^{-1}(Ts^1_0 (T\Delta^0)) \to A$ can
be written as $a(t,s)dt + b(t,s)ds$ over the base map $\gamma(t,s)$
after a suitable choice of parametrization\footnote{We need the one that 
$\gamma(0,s)=x$ and $\gamma(1,s)=y$ being constant for all $s\in[0,1]$. } of
the disk $(d^2_2)^{-1} (s^1_0( \Delta^0))$.
Then we naturally have $b(0,s)=f(0,s)(\frac{\partial}{\partial s})=0$ and
$b(1,s)=f(1,s)(\frac{\partial}{\partial s})=0$.  Moreover the morphism $f$ is a
Lie algebroid morphism if and only if $a(t,s)$ and $b(t,s)$ satisfy 
equation  \eqref{eq:hom} which defines the $A$-homotopy (see for example \cite[Cor. 4.3]{brahic-zhu}). 
\end{itemize}
\begin{remark} We do not know whether $S(A)$ is a simplicial
manifold or further a {\em Kan simplicial manifold}, namely a Lie
$\infty$-groupoid as in Definition \ref{def:defngroupoid},  though it is so for
a Lie algebra \cite{henriques}. Unlike that $S_1(A)$ being a Banach manifold involves solving an ODE, it is
not clear how to solve directly the corresponding PDE for $S_2(A)$ to
be a Banach manifold. This is one of the open questions left at the
end of \cite{henriques}. 
\end{remark}

However, although $S(A)$ has a clear geometric meaning,
it involves infinite dimensional manifold, and although $\cG(A)$ is an
\'etale SLie groupoid, the 2-truncation $\tau_2(S(A))$ will not
be  2-\'etale\footnote{See Definition \ref{def:2etale}.}.   
To achieve the
\'etale version we need to use a sub-simplicial set based on the
good \'etale covering $G_0(A)$ of $\cG(A)$ in Lemma \ref{lemma:etale-good}, namely
\begin{equation}\label{eq:etale-simp} 
S^{\acute{e}}_0(A)=M, \;  S^{\acute{e}}_1(A)=G_0(A),  \;  S^{\acute{e}}_i(A)= S_i (A)|_{G_0(A)},
\; \text{ for} \; i\geq 2, \end{equation}
where 
\[ S^{\acute{e}}_i (A)|_{G_0(A)}:= \{ x \in S_i(A): \text{$1$-skeleton of
  $x$ is made up by elements in $G_0(A)$.} \} \]

\subsection{Lie 2-groupoid version}\label{sec:lie2-v}
Now we look for the Lie 2-groupoid corresponding to $\cG(A) \rra M$. For
this we need the notion of a good chart.  

\begin{defi}
If a map $\bar{f}$ from a manifold $N$ to
an \'etale
differentiable stack $\cX$ with  an \'etale chart $X_0$ comes from an embedding $f: N\to X_0$, we call $X_0$ a
{\em good \'etale chart} and its corresponding groupoid presentation
$G:= G_1 \rra G_0$ a {\em good \'etale groupoid presentation}. 
\end{defi}

Then we have the following technical lemma (proved in \cite[Section~3]{z:tgpd-2}),

\begin{lemma}\label{lemma:embedding}
Given an immersion $\bar{e}: M\to \cG$ from a manifold $M$ to an
\'etale stack $\cG$,  there is a good \'etale chart $G_0$ of $\cG$,
i.e. $\bar{e}$ can be lifted to an embedding $e: M \to G_0$.  
\end{lemma}

The \'etale chart $G_0(A)=\sqcup_i P_i$ of $\cG(A)$, made
up by local transversals $P_i$ of the foliation $\cF$ on $P_0 A$, is
usually not a good \'etale chart itself. But we have the following lemma:

\begin{lemma}  \label{lemma:etale-good}
The \'etale chart $G_0(A)= \sqcup_i P_i$ of $\cG(A)$ can be completed to a
good \'etale chart of $\cG(A)$.
\end{lemma}
\begin{proof}
We recall the construction of the
local groupoid $G^{loc}(A)$ of $A$ in \cite{cf}. Take a small tubular
neighborhood\footnote{By \cite[Chap.VII.4]{lang}, a tubular neighborhood can be taken isomorphic to a vector bundle over the  submanifold if the submanifold is paracompact  and the ambient manifold is a Hilbert manifold, which is our case (see Remark \ref{rk:hilbert}).} $\cO$ of $M$ in the $A_0$-path space $P_0 A$ so that
the foliation $\cF$ restricted to $\cO$, which we denote by
$\cF|_{\cO}$, has good transversal sections, that is all the leaves in
$\cF|_\cO$ intersect every transversal section only once. Take $\cO_1 \subset \cO$ a smaller tubular neighborhood of $M$ such that $\cO_1 \times_M \cO_1 \xrightarrow{\odot} \cO$, and for all $a(t)\in \cO_1$, the inverse path $a(1-t)$ is also  in $\cO_1$.  Then the
quotients  $U:=\cO/(\cF|_\cO)$ and $V:=\cO_1/(\cF|_\cO)$ yield a local Lie
groupoid over $M$.

The open set $\cO/(\cF|_\cO)=:U$ can be also visualized by gluing $P_i
|_\cO$ together via the equivalence induced by $\cF|_{\cO}$ on
$P_i|_{\cO}$. Then, by definition, there is a natural embedding $U \to \sqcup_i P_i|_{\cO} \to \sqcup_i P_i$. Composing with the representable submersion $\sqcup_i P_i \to \cG$, we obtain a representable submersion $U\to \cG$. 

Although $G_0(A)$ is not necessarily a good chart, we join $U$ to
$G_0(A)$, that is, we redefine $G_0(A):= U
\sqcup (\sqcup_i P_i)$. Then $G_0(A)$ is still a chart of $\cG(A)$ since $U\to \cG$ is a representable submersion.  Thus $G_0(A)$ becomes good and the \'etale groupoid
$G_1(A):=G_0(A)\times_{\cG(A)} G_0(A) 
\rra G_0(A)$ via $G_0(A) \to \cG(A)$ is a good \'etale groupoid presentation. 
\end{proof}

To avoid redundant notation, from now on $G_0(A)$ denotes this completed good \'etale
chart of $\cG(A)$ and $G(A)$ this good \'etale groupoid presentation. Then the multiplication bimodule $E_m^{G(A)}$ of $G(A)$
is 
\[ E_m^{G(A)}=\big( G_0(A) \times_M G_0(A)  \big) \times_{\odot, P_0
  A, \bt }
Mon (P_0 A)_1 \times_{\bs, P_0 A, i } G_0(A), \]
 with $G_0(A) \times_M G_0(A) \xrightarrow{\odot} P_0 A$ given by $(a, b) \mapsto a \odot b$
as in \eqref{eq:con-2} and the immersion $G_0(A) \xrightarrow{i}
P_0 A$.  Since $Mon(P_0 A)_1 \times_{P_0 A} G_0(A) \to P_0 A$ is surjective and  \'etale (see \cite[Example 5.9]{moerdijk}), we verify
directly that $E_m^{G(A)}$ is a manifold. Thus by the 1-1
correspondence in \eqref{eq:s-2}, 
\begin{equation}\label{eq:2gpd-ft}\xymatrix{\Gamma(A):  E_m^{G(A)} \ar@<1ex>[r]  \ar@<-1ex>[r]  \ar[r]  & G_0(A)
  \ar@<1ex>[r]  \ar[r]  & M} \end{equation}
is a Lie 2-groupoid. We call $\Gamma(A)$ the universal Lie 2-groupoid integrating $A$. Notice that different choices of $G_0(A)$ give the same Lie 2-groupoid up to 1-Morita equivalence (see \cite[Prop.-Def.2.13]{z:tgpd-2}.

There is also a
natural infinite dimensional Lie 2-groupoid corresponding to $\cG(A)
\rra M$ obtained by using the chart $P_0 A$, 
\begin{equation}\label{eq:inf-dim} \xymatrix{ Y(A): E_m \ar@<1ex>[r]  \ar@<-1ex>[r]  \ar[r]  & P_0 A
  \ar@<1ex>[r]  \ar[r]  & M}, \end{equation}
where $E_m=(P_0 A \times_M P_0 A) \times_{\odot, P_0 A, \bt}
Mon_{\cF}(P_0 A)$ is the multiplication bimodule of $\cG(A)$ with
respect to the chart $P_0
A$.  In fact \eqref{eq:2gpd-ft} is simply the
pullback 2-groupoid of \eqref{eq:inf-dim} via the
map $G_0(A) \xrightarrow{i} P_0A$. Even though $G_0(A) \xrightarrow{i} P_0A$ is not a surjective submersion, these two 2-groupoids are Morita equivalent because both of them correspond to the same stacky groupoid $\cG(A) \rra M$. 

Then Theorem \ref{thm:2-a} is equivalent to the following:
\begin{thm}\label{thm:2-a'}
The 2-truncation $\tau_2(S^{\acute{e}}(A))$ is exactly the 2-\'etale\footnote{see Definition \ref{def:2etale}} Lie 2-groupoid $\Gamma(A)$ defined in \eqref{eq:2gpd-ft}, and the 2-truncation $\tau_2(S(A))$ is exactly the Lie 2-groupoid $Y(A)$ defined in \eqref{eq:inf-dim}.  Hence both $\tau_2(S^{\acute{e}}(A))$ and  $\tau_2(S(A))$  corresponds to $\cG(A)$ with the
correspondence in \cite{z:tgpd-2}.
\end{thm}
\begin{proof}
To prove $\tau_2(S^{\acute{e}}(A))\cong \Gamma (A)$, we only need to show that $\tau_2(S^{\acute{e}}(A))_2 \cong E^{G(A)}_m$ and the rest is
easy to check. Thanks to the boundary condition in the definition
of $P_0 A$,  $Mon(P_0 A)_1$ is a quotient of the space of
bigons in $S_2(A)$.   Thus we have, 
\begin{equation}\label{eq:em}
  E^{G(A)}_m = S^{\acute{e}}_2(A)/\sim,  
\end{equation}
where  $a\sim b$ in $S^{\acute{e}}_2(A)$ if and only if there is a
path $a(\epsilon)$ in $S^{\acute{e}}_2(A) $ with fixed boundary,
i.e. $\partial(a(\epsilon))=\partial a = \partial b$,   such that $a(0)=a$ and
$a(1)=b$.   By   \cite[Lemma 4.5]{brahic-zhu}, a path of
$A$-homotopies  with fixed boundary can be connected to make a
Lie algebroid morphism $T\Delta^3 \to A$ in a unique way.
Thus  what we mod out in \eqref{eq:em} is
exactly an element in $S_3(A)$. Hence   $\tau_2(S^{\acute{e}}(A))
\cong E^{G(A)}_m$. The proof of the statement $\tau_2(S(A))\cong Y(A)$ is similar.
\end{proof}

\section{Lie II theorem}
In this section we prove Theorem \ref{thm:lieii}. The major tool is
the adaptation of Quillen's small object argument for simplicial
manifolds. It is necessary to demonstrate this since simplicial manifolds do not form a
model category. 

\subsection{Universal Lie 2-groupoids via local groupoids}
Although there is no one-to-one correspondence of Lie algebroids and
Lie groupoids, Lie algebroids and local Lie groupoids do have  a one-to-one
correspondence. Given a Lie algebroid $A$, its corresponding local Lie
groupoid $G^{loc}(A)$ is constructed explicitly as in Lemma
\ref{lemma:etale-good}. 
Then we naturally have a 2-groupoid
$\tau_2(Kan(NG^{loc}(A)))$, and we prove:

\begin{prop}\label{prop:glocalg}
Given a Lie algebroid $A$,  $\tau_2(Kan(NG^{loc}(A)))$ is a Lie 2-groupoid and is Morita equivalent to the universal Lie
 2-groupoid $\Gamma(A)$ associated to $A$. 
\end{prop}

We prove this proposition via infinite dimensional spaces because there we
can to a certain extent avoid the problem of the ``non-existing strict
morphism'' (see \cite[Section 3.2]{z:kan}) and simplify considerably the proof.

We use the same notation as in Lemma \ref{lemma:etale-good}, and we take $X(A)$ to be the pull-back simplicial manifold of $NG^{loc}(A)$
along the surjective submersion $\cO_1 \to V$. Then by construction, $X(A)$ satisfies properties
A and B. Since $X(A)_1$ contains small paths and $X(A)_2$ contains
small triangles,  $X(A)_n$ is a neighborhood of the identities
$Y_0(A)=M$ in $Y(A)_n$. Thus $X(A)$ can be viewed as a local part of $Y(A)$.

\begin{lemma} \label{lemma:local}
The 2-truncation $\tau_2(Kan(X(A)))$ is a Lie 2-groupoid and is Morita equivalent to $Y(A)$.  
\end{lemma}
\begin{proof}
To simplify the notation, we omit ``$(A)$'' in the proof of this
Lemma. 

We first construct a morphism $f': Kan(Y) \to Y  $ inductively thanks to
\eqref{eq:iso-colim}.   The construction is similar to the
construction in \cite[Theorem 3.6]{z:kan}. There we could not construct
a global map. However here we are able to do so because we can choose a
global section of $\hom(\L[2, 1], Y) \to Y_2$, for example the one given by the ``flat''
triangle of concatenation ~\eqref{eq:con-2}, 
\begin{equation}\label{eq:flat}
P_0 A \times_M P_0 A \to
\hom(T\Delta_2, A), \quad (a, b) \mapsto \text{the flat triangle whose
  three sides are}\; a, b, a\odot b. \end{equation}

Clearly there is an embedding $ X
\xrightarrow{\iota} Y$, hence a map $Kan(X) \to Kan(Y)$.  Thus we
have a morphism
$f: Kan(X) \to Y$. 

We do not have a hypercover directly from  $\tau_2(Kan(X))$ to
$Y$ because $Kan(X)_1 \to Y_1$ is not surjective. Thus we will use
Lemma~\ref{lemma:me} and verify that the conditions in Remark \ref{rk:needed} are satisfied.

Given a triangle, we add a point and divide it into three
triangles. Among these three triangles, one can be a weird one in $X_1\times_{X_0}X_1$, as
\eqref{pic:x12} shows. A triangulation is {\em good} if it is made of a sequence
of such divisions.  Thus, more geometrically,  an element in
$Kan(X)_2$ is a set of  small
triangles of $X_2$ matching together in a good triangulation
situation.  Then $f_1: Kan(X)_1 \to Y_1$ operates by concatenating 
small paths into a long path and $f_2$ by  composing small triangles into a big one in $Y_2$ with the flat
filling of $X_1\times_{X_0} X_1$ as in \eqref{eq:flat}. The embedding $f_1|_{X_1}: X_1 \to Y_1$ is a submersion in particular.

Recall that we need an intermediate $Z$ with $Z_0=M$ and 
\[Z_1  = \\hom(\Lambda[2,1], Kan(X)) \times_{\hom(\Lambda[2,1], Y)} Y_2\sqcup X_1\times_{f_1, Y_1, d_0} b(Y_2) \]
The map $h'_1: \hom(\Lambda[2,1], Kan(X)) \times_{\hom(\Lambda[2,1], Y)} Y_2
\to Y_1$ is surjective because any $A_0$-path $a_1 \in Y_1$ is
$A$-homotopic to a concatenation of small paths in $X_1$. More
precisely, there are elements $a_2, a_0 \in Kan(X)_1$ such that
$f_1(a_0) \odot f_1(a_2)$ is $A$-homotopic to
$a_1$. Then this $A$-homotopy gives an element $\eta\in Y_2$ with
three sides $f_1(a_0)$, $a_1$, and   $f_1(a_2)$
respectively. Thus $h'_1(a_0, a_2, \eta) = a_1$. Notice that given a
small path $a_1^s$ around $a_1$,  the choice
of $(a^s_0, a^s_2, \eta^s)$ can be made smoothly depending on $s$. Thus
$h'_1$ is also a submersion.

The morphism $\mu: \hom(\Lambda[2,1], Kan(X)) \to Kan(X)_1$ is obtained by taking the limit of the natural embedding $X^\beta_1 \times_{X_0} X^\beta_1 \to X^{\beta+1}_1$. It is clear that $\mu$ is a submersion, thus the composed map $g'_1$ is a submersion since the Kan projection $Y_2 \to \hom(\Lambda[2,1], Y)$ is a submersion. Then joining with $pr_{X_1}$,  the  map $g_1=g'_1\sqcup pr_{X_1} :  Z_1 \to Kan(X)_1$ is a surjective submersion. 

Now we need to verify that the natural map obtained by composing small
triangles via the Kan condition
\[ Kan(X)_2 \xrightarrow{\phi} \hom(\partial \Delta [2], Kan(X)) \times_{f, \hom(\partial \Delta[2], Y)} Y_2, \]
induces an isomorphism when passing to the quotient of the 2-truncation. Then $\tau_2(Kan(X))$ is a Lie 2-groupoid and Morita equivalent to $Y$.

Thanks to the weird triangle in \eqref{pic:x12}, we have a sort of
division as in the second group in $\hom(\L[3, j], X^1)$ in
\eqref{pic:x22}. After several times,
we can always end up with small triangles with at most $\frac{1}{2}$
perimeter of the original one. Hence any  $y \in Y_2$ can be divided
into small triangles that lie in $X_2$.   
Thus the map $\phi$ is surjective up to the homotopy $\sim_2$ given by elements in $Y_3$.  

On the other hand, suppose two relatively homotopic elements $(\partial z, y)$, $(\partial z', y')$ $\in$ $\hom(\partial \Delta [2], Kan(X))$ $\times_{\hom(\partial \Delta [2], Y)}$ $ Y_2 $ have preimages $x, x'$, i.e. two
triangulations $x, x'$ of two homotopic triangles $\eta, \eta' \in
\hom(T\Delta_2 , A)$  sharing the same division of the border. Since $\eta$ and
$\eta'$ can be connected by a path $\eta_t \in \hom(T\Delta_2, A)$, we
might assume that this path is a small path and we
can deform $x'$ to make it also a triangulation of $\eta$. Then  $x$ and $x'$ have a
common
subtriangulation $\tilde{x}$ and $\tilde{x}'$ (though with a possibly
different order of composition) which has the same division of the border.  However different order of
composition differs by an element in $Kan(X)_3$. A triangulation
and its subtriangulation sharing the same border
also differ by an element in $Kan(X)_3$. Hence $\phi$ is injective up to the homotopy $\sim_2$ given by elements in $Kan(X)_3$.  
\end{proof}

\begin{proof}[Proof of Prop. \ref{prop:glocalg}]
We take the Lie 2-groupoid $Y(A)$  and $X(A)$ as above.  Since $X(A) $ is the pull-back of $N G^{loc}(A)  $ by the surjective submersion $\cO_1 \to V$, the map $X(A) \to NG^{loc}(A)$
is a hypercover.  Hence
$\tau_2(Kan(X(A)) \to \tau_2(Kan(NG^{loc}(A)))$ is a hypercover by 
\cite[Thm. 3.11]{z:kan}. Thus we have a composed Morita
equivalence, 
\[ \Gamma(A) \me Y(A) \stackrel{\sim}{\leftrightarrow} \tau_2(Kan(X(A)))
\stackrel{\sim}{\to} \tau_2( Kan (NG^{loc}(A))) ,  \] by Remark \ref{rk:finite} and Lemma \ref{lemma:local}.
\end{proof}

\subsection{Differentiating back to Lie algebroids}

\begin{defi}\label{def:2etale}
As in \cite{z:tgpd-2}, we call a Lie 2-groupoid $X$ {\em 2-\'etale} if  the Kan projections $X_2 \to \hom(\Lambda[2,j], X)$ is \'etale for  $j=0,1,2$.
\end{defi}

In a 2-\'etale Lie 2-groupoid,  $X_0$ is an embedding submanifold of $X_1$ and $X_2$ via compositions of
degeneracy maps. Moreover $Kan(2,1)$ tells us that $X_2 \to X_1 \times_{d_0, X_0, d_2}
  X_1 $ is \'etale. Hence there is a tubular neighborhood
$V \subset X_1$ of $s_0: X_0\to X_1$ such that there is a unique section $\sigma_m: V
\times_{d_0, X_0, d_2} V \to X_2$ extending the section given by the
degeneracy map  $X_0 \times_{d_0,X_0, d_2} X_0 \to X_2 $.  The image $d_1 \circ
\sigma_m( V\times_{d_0, X_0, d_2} X_0 )$ equals  $V$, hence if we take
  $V$ small enough, we have another open neighborhood $U$ of
  $X_0$ in $X_1$ satisfying  $V  \subset  d_1 \circ
  \sigma_m ( V
\times_{d_0, X_0, d_2} V  )\subset U$. We define
\[m=d_1\circ \sigma_m: V\times_{d_0, X_0, d_2} V \to U.  \]

The identity embedding $X_0 \to U$ is given by $s_0$. 

Finally, $Kan(2,0)$ tells us that $X_2 \to X_1\times_{d_2, X_0, d_1} X_1 $
is \'etale. Hence there is a unique section $\sigma_i: U\times_{d_2,
  X_0, d_1} X_0 \to X_2$ extending $X_0\times_{d_2,
  X_0, d_1} X_0 \to X_2$. We take $U$ small
enough (for example $U\cap d_1\circ \sigma_i(U \times_{X_0} X_0)$), and define the composed map
\[i: U \xrightarrow{ id \times d_1} U \times_{d_2,
  X_0, d_1} X_0 \to X_2 \xrightarrow{\sigma_i} U \]
to be the inverse map. Then $m(i(g), g)=1$ and  $m(g, i(g))=1$ by construction.

Then the uniqueness of the sections $\sigma_m$ and $\sigma_i$ gives us
strict Kan condition $Kan(2, j)!$ locally. Hence $V \subset U
\rra X_0$
with the above structure maps is a local Lie groupoid.  The
structure of the local groupoid does not depend on the choice of
$U$ since two local groupoids are defined  to be isomorphic if they are the same in
an open neighborhood of $X_0$ in $X_1$. We call it the {\em  local Lie
groupoid associated to the 2-\'etale Lie 2-groupoid $X$} and denote it by $X^{loc}$. By construction $(NX^{loc})_n$ is a tubular neighborhood of the identity $X_0$ in $X_n$. Thus we naturally have an inclusion $NX^{loc} \subset X$ of simplicial manifolds.

Similarly given a W-groupoid $\cG \rra M$ where
$\cG$ is presented by $G_1 \rra G_0$, there are neighborhoods
$V\subset U \subset G_0$ of $M$ such that all the groupoid structure
maps descend to $(U, V)$ and make $(U, V)\rra M$ a local Lie groupoid 
\cite[Section~5]{tz}\footnote{Thanks to Lemma
\ref{lemma:embedding} we can construct the local groupoid at once
and it is not necessary to divide $M$ into pieces as it is done in the
cited reference.}. We
call it  the {\em  local Lie
groupoid associated to W-groupoid $\cG \rra M$} and denote it by $G^{loc}$.

The Lie algebroid of $X^{loc}$ (respectively $G^{loc}$) is defined
to be {\em the Lie algebroid of the 2-\'etale Lie 2-groupoid} $X$
(respectively {\em W-groupoid} $\cG$).

\begin{ep}
If we take the universal 2-\'etale Lie 2-groupoid $\Gamma(A)$, the local
Lie groupoid associated to it is exactly $G^{loc}(A)$.  
\end{ep}

A {\em W-groupoid morphism} $\Phi: (\cG\rra M) \to (\cH\rra N)$ is
made up by a map $\Phi_1: \cG\to \cH$ between stacks and by a map $\Phi_0:
M\to N$ such that they  preserve the W-groupoid structure maps up
to 2-morphisms, which  should satisfy again higher
coherence conditions linking the 2-commutative diagrams of $\cG$ and $\cH$
(also see \cite[Section~4]{z:tgpd-2}). Such a morphism, in the world of 2-\'etale Lie 2-groupoids, corresponds to a generalized
morphism $\Psi: \gm{X}{Y}$ via $Z$ such that the morphism $Z\to X$ makes $Z_0
\cong X_0$. We call such generalized morphisms {\em special generalized morphisms}.  Now we differentiate this generalized morphism to a
Lie algebroid morphism.

First we notice that $Z_1 \to X_1$ is a surjective submersion. Thus we
choose a local section $\sigma_1$ of $Z_1 \to X_1$ from a small neighborhood of $X_0$ in
$X_1$ extending the isomorphism $X_0 \to Z_0$.  Then we obtain a local
morphism $\sigma_2: X^{loc}_2 \to Z_2$ by $x_2 \mapsto (x_2, \sigma_1(d_0(x_2)),
\sigma_1(d_1(x_2)), \sigma_1(d_2(x_2)))$ since $Z_2\cong
X_2\times_{\hom(\partial \Delta[2], X)} \hom(\partial \Delta[2],
Z)$. Similarly we obtain higher $\sigma_i: X^{loc}_i \to Z_i$. It is
easy to see that with $\sigma_0: X_0 \to Z_0$ as the isomorphism,
$\sigma: NX^{loc} \to NZ^{loc}$ is a simplicial morphism. Thus the composed
simplicial morphism $NX^{loc} \to NZ^{loc} \to NY^{loc} $ gives a
local group morphism $\phi: X^{loc} \to Y^{loc}$ since the local
groupoid structure is determined by the simplicial manifold structure
of the nerve. Then differentiating $\phi$ gives us a Lie algebroid
morphism $A \to B$ where $A$ and $B$ are the Lie algebroid of $X$ and
$Y$ respectively. We call $\Psi$, or the corresponding W-groupoid
morphism $\Phi$, the
{\em integration of} $\phi$.

\begin{lemma} \label{lemma:lieii}
If $\phi$ is a Lie algebroid morphism $A \to B$, then it induces a
Lie 2-groupoid morphism $\Phi: \Gamma(A) \to \Gamma(B)$ integrating
the Lie algebroid morphism $\phi$. 
\end{lemma}
\begin{proof}
A  Lie algebroid morphism $\phi: A\to B$ induces a simplicial
morphism $\phi_i: S(A)_i \to S(B)_i$. Thus it gives a Lie 2-groupoid
morphism on the 2-truncations. By Theorem \ref{thm:2-a'}, we obtain a
Lie 2-groupoid morphism $\Phi: \Gamma(A) \to \Gamma(B)$. Then by what we discuss
above, $\Phi$ gives a local
groupoid morphism $\Phi^{loc}: G^{loc}(A)\to G^{loc}(B)$ which
maps equivalence classes of $A$-paths in $A$ to those in $B$ via
$\phi$. Note that the local groupoids can be understood as equivalence
classes of $A$-paths as explained in Lemma \ref{lemma:local}. Therefore the corresponding Lie
algebroid map is exactly $\phi$. 
\end{proof}

By the one-to-one correspondence between the W-groupoid and 2-\'etale
Lie 2-groupoid, Theorem \ref{thm:lieii} is equivalent to the following theorem:
\begin{thm}
Given a Lie algebroid morphism $A \xrightarrow{\phi} B$, let $G(A)$ be the universal
Lie 2-groupoid of $A$, and $W$ any 2-\'etale Lie 2-groupoid
corresponding to $B$.
Then there is a special generalized morphism $\Psi: \xymatrix@C=.45cm{G(A)
  \ar@{.>}[r] & W}$ integrating $\phi$. This generalized morphism is unique up to
2-morphisms. 
\end{thm}
\begin{proof}
To build the map $\Psi$,
by Lemma \ref{lemma:lieii}, we only have to treat the situation
when $\phi=id: A\to A$, that is, if a Lie 2-groupoid $W$ has
algebroid $A$, then there is a generalized morphism $\Psi: \Gamma(A) \to
W$ integrating $id: A\to A$.  Then $W^{loc}$ is the local Lie
groupoid of $A$ since the Lie algebroid corresponding to $W$ is $A$.
Thus we have an embedding of simplicial manifolds
$NG^{loc}(A)=NW^{loc}\subset W$. By Prop. \ref{prop:glocalg} and Thm. \ref{lemma:kankan},  we obtain a composed generalized morphism
\[ \Psi: \Gamma(A) \stackrel{\sim}{\leftrightarrow} \tau_2( Kan (NG^{loc}(A))) \to \tau_2(Kan(W)) \stackrel{\sim}{\leftrightarrow} \tau_2(W) =W.  \]
Locally $\Psi$ is the
identification map $G^{loc}(A) \to W^{loc}$, thus it integrates $id:
A\to A$.  Different choices of $G^{loc}(A)$ give us the same generalized morphism up to a  2-morphism.

On the other hand, if we have two special generalized morphisms
$\Psi, \Psi': \xymatrix@C=.45cm{\Gamma(A) \ar@{.>}[r]& H}$ integrating the Lie algebroid morphism
$\phi: A \to B$, via $Z$ and $Z'$, then locally we can
invert the projections $Z \to \Gamma(A)$ and $Z' \to \Gamma(A)$ and get two strict
morphisms $\Psi^{loc}, \Psi'^{loc}: NG^{loc}(A) \to H$. By choosing
$G^{loc}(A)$ close enough to the identity section,
$\Psi^{loc}=\Psi'^{loc}$ because they both integrate the Lie algebroid morphism 
$\phi: A\to B$.  Then since the local
morphism determines the global generalized morphism,
\[ \Gamma(A) \stackrel{\sim}{\leftrightarrow}  \tau_2(Kan(G^{loc}(A)))
\xrightarrow{\tau_2(Kan(\Psi^{loc}))=\tau_2(Kan({\Psi'}^{loc}))} 
  \tau_2(Kan(H)) \stackrel{\sim}{\leftrightarrow}  \tau_2(H)=H,\] 
so we must have $\Psi =\Psi'$ up to a 2-morphism. 
\end{proof}

\section{Connectedness}\label{sec:conn}
In this section, we will show that $Y(A)$ and $\cG(A)$ are source-2-connected,
that is, their source fibres have trivial homotopy groups  $\pi_{0}$, $\pi_1$ and $\pi_2$.  The idea is that
the $n$-truncation of $S(A)$ is the universal object integrating $A$ in the world
of differentiable $(n-1)$-stacks or equivalently Lie $(n-1)$-groupoids. On the other hand, if we take a
$k$-connected cover of this $(n-1)$-stack for $k=1,
..., n$, the cover is also a universal object. Hence we must have: 

{\em  The universal object in differentiable
$(n-1)$-stacks has to be $n$-connected, by the same argument according
to which the
universal group for a Lie algebra is $1$-connected.}

But to make this into a proper proof, to construct this $k$-connected cover, we must be equipped with  homotopy
theory for Lie $n$-groupoids or differentiable $n$-stacks. Notice that even for a manifold $M$, its $k$-connected cover
might be a higher stack (see Section \ref{sec:ep}).  Since in this article we are aiming at  2-groupoids, we give
a direct proof of the 2-connectedness using the above idea. We also give some examples of 2-connected covers in the
example section. In this method, we view the universal integrating
object as a (1)-stack with an additional groupoid structure. Thus it
is more natural to use the stacky Lie groupoid model for the universal
object.

Homotopy groups of simplicial schemes \cite{friedlander},
simplicial sets \cite{may}, orbifolds \cite{bri-hae}, and
topological stacks \cite{noohi:top} are well-known or have been
recently studied. They have properties that parallel those of the
usual homotopy groups of topological spaces; for
example, $\pi_n$ is a group. Moreover, we believe that these
theories are all the same under suitable equivalences, but
sometimes not the obvious one (see \cite{friedlander}). For
example, in the sense of \cite{may}, the simplicial homotopy
groups $\pi^s_n$ of the $2$-truncation of $S(A)$ is automatically
trivial by definition when $n\geq 2$. But we will show in Remark \ref{rk:2-conn-cover},
$\pi_2 (\cG(TS^2))\neq 0$ with $\pi_2$ defined in Definition \ref{def:pi}. Indeed they are two {\em different}
though related types of homotopy groups if applied to simplicial
manifolds.  Here we give a version of
homotopy groups of differentiable stacks which generalizes that of
orbifolds in \cite{bri-hae}.

As before we denote 2-morphisms as $\sim$.
\begin{defi} \label{def:pi}At a point $x: pt \to \cX$, the $n$-th
homotopy group is defined as, 
\[\pi_n(\cX, x) = \{ f: S^n\to \cX,\,\text{ such
that}\, f|_N \sim i_x \}/ \text{homotopies}, \] here $i_x$ is the
morphism from the north pole $N= pt \overset{x}{\to} \cX$, and
$f_0$ is homotopic to $f_1$ if and only if there is $F: S^n\times [0,1] \to
\cX$ such that
\[ F(N, s)\sim i_x, F(\cdot , 0) \sim f_0, F(\cdot , 1)\sim f_1.\]
\end{defi}

Since in the main results of this paper, we are using nothing more than the
definition of the homotopy groups, we omit a detailed
discussion of homotopy groups.

Source connectedness for $\cG(A)$,
i.e. $\pi_0(\bs^{-1}(x), x)=0$ is easy to see. It is implied by the
fact that every $A_0$ path $a(t)$ can be connected to a constant path
$0_x$ via $sa(st)$ with $s\in [0,1]$. Here $x=a(0)$. Now we set off to
prove the less trivial part about $\pi_1$ and $\pi_2$.

We recall how we prove that the universal Lie group $G(\g)=P_0 \g/\sim$ of a Lie algebra
$\g$, constructed by $A$-paths modding out $A$-homotopies in $\g$, is
simply connected. We use only Lie II theorem
but {\em without using} a covering space.

Given a loop $S^1 \xrightarrow{\varphi} G$, we have a groupoid
morphism $P(S^1) \xrightarrow{\Phi_{gpd}} G$  given by $\Phi_{gpd}
(c, d) = \varphi(c) \varphi(d)^{-1}$. Here $P(\square)$ denotes the
pair groupoid $\square \times \square \rra \square$. We use this
notation from now on. The corresponding Lie algebroid
morphism $TS^1 \xrightarrow{\phi_{algd}} \g$ is in fact a
$\g$-loop $\phi_{algd}: \varphi'(t)\varphi(t)^{-1}$. By re-parametrization, we might as well assume
that $\dot{\varphi}([t])=0$, when $t=0,1$. Here we view $S^1=\R/\Z$
and $[\cdot]$ denotes an equivalence class. Then
$\phi_{algd} $ is in $P_0\g$, and $\varphi$ is
contractible if and only if $\phi_{algd}$ is contractible,
i.e. $\phi_{algd}\sim 0$ in the quotient $G(\g)=P_0\g/\sim$. In fact, if
$\phi_{algd} \sim 0$ via a $\g$-homotopy $F_{\phi}: TD^2 \to \g$, we
integrate it to a homotopy $F_{\Phi}:P( D^2)  \to G$. Then
$F_{\Phi}$ restricted to one copy of $D^2$ gives the homotopy which
contracts $\varphi$. 

To prove $\phi_{algd}\sim 0$, we use the fact that $\phi_{algd}$ is a 
$\g$-path coming from the groupoid
morphism from the pair groupoid $P(S^1) \xrightarrow{\Phi_{gpd}} G$, with $\Phi_{gpd}
(c, d) = \varphi(c) \varphi(d)^{-1}$. More precisely, the
integration morphism $ \tilde{\Phi}_{gpd}$ of $\phi_{algd}$ factors
through $P(S^1)$,  
\[ \tilde{\Phi}_{gpd}: G(TS^1) = (\R \times \R /\Z \rra S^1) \xrightarrow{pr}
P(S^1)\xrightarrow{\Phi_{gpd}} G,  \]
with every arrow a strict groupoid morphism. This is because the infinitesimal morphism of both $ \tilde{\Phi}_{gpd}$ and
$\Phi_{gpd}$ is $\phi_{algd}$. (Here we use Lie II theorem). 

Taking $[(0,1)]$ and $[(0,0)]$ in  $\R \times \R /\Z$,  we have 
\begin{equation}\label{eq:01-00} \tilde{\Phi}_{gpd}([(0,1)])= \tilde{\Phi}_{gpd}([(0,0)]), \end{equation}
because $pr([(0,0)])=pr([(0,1)])=([0],[0]) \in S^1\times S^1 = \R/\Z \times
\R/\Z$. Here $[\cdot]$ denotes equivalence classes in various
quotients. We define maps $\gamma_{0,1,2}: S^1\to S^1$ by
$\gamma_0([t])=[0]$, $\gamma_1([t])=[\tau(t)]$, and
$\gamma_2([t])=[t]$, with $t\in [0,1]$ and 
\begin{equation}\label{eq:tau}
\tau:[0,1]\to
[0,1],\quad \text{ which satisfies
$\tau'(0)=\tau'(1)=0, \tau'(t)\geq 0 $. } 
\end{equation} 
As $A$-paths, $T\gamma_1, T\gamma_2: TS^1 \to TS^1$ differ by a
re-parametrization $\tau$. In particular,  $T\gamma_{1}$ and $T\gamma_2$ are
$A$-homotopic. Thus  $T\gamma_1$ represents $[(0, 1)]$ and $T\gamma_0$ represents $[(0,0)]$
in $\cG(TS^1)$ respectively.  
  Then \eqref{eq:01-00} tells us that, on the algebroid level, $\phi_{algd}
\circ T\gamma_1 \sim \phi_{algd}\circ T\gamma_0$ are two equivalent
$\g$-paths.  

But it is clear that $\phi_{algd}\circ T\gamma_1 = \phi_{algd}$ and
$\phi\circ T\gamma_0 =0$ the $0$-path. Hence $\phi_{algd} \sim 0$. 

Now we use the same idea to prove our theorem. 

\begin{proof} [Proof of Thmeorem \ref{thm:pi2}]
It is almost the same proof to show that  $\cG(A)$
has simply connected source fibre $\bs^{-1}(x)$ for all $x\in M$. We
only have to replace $\g$ by $A$ and insert the base point $x$. We omit the
proof here. 

We focus on adapting the proof for the
2-connectedness. Given $x\in M$, let
$\varphi$ be a smooth map $(S^2, N) \to (\bs^{-1}(x), x)$.  Here $x:
pt \to \bs^{-1}(x)$ denotes also the base point. Then  similarly, we form a
W-groupoid morphism $\Phi_{gpd}:P(S^2) \to \cG(A)$ by $\Phi_{gpd}(c,
d)=\varphi(c)\varphi^{-1}(d)$ with the base map
$\Phi^0_{gpd}(c)=\bt(\varphi(c))$. Then $\Phi_{gpd}$ induces a Lie algebroid
morphism $\phi_{algd}: TS^2 \to A$.

For any Lie algebroid $A$, there is another stacky Lie groupoid $\cH(A):=\cH(A)_1 \rra M$ with arrow space $\cH(A)_1$ presented by the holonomy groupoid $Hol(P_0 A)\rra
P_0 A$ of the same foliation. By \cite[Theorem 1.3]{tz}, when a Lie algebroid $A$ is integrable to a source simply-connected Lie groupoid $G$, $\cH(A)$ is simply $G$. Since the source simply-connected Lie groupoid of  $TS^2$ is $P(S^2)$, the stacky Lie groupoid $\cH(TS^2)$ is this Lie groupoid $P(S^2)$. 

Then the projection $Mon(P_0 TS^2) \to Hol(P_0 TS^2)$ induces a stacky
groupoid morphism $pr: \cG(TS^2) \to P(S^2)$ whose corresponding
infinitesimal morphism is $id: TS^2 \to TS^2$. The composed stacky
groupoid morphism \[ \tilde{\Phi}_{gpd}: \cG(TS^2)\xrightarrow{pr} P(
S^2) \xrightarrow{\Phi_{gpd}} \cG(A),  \] still corresponds
infinitesimally to $\phi_{algd}$. Thus $\phi_{algd}$ integrates to a W-groupoid morphism $\tilde{\Phi}_{gpd}: \cG(TS^2)\to \cG(A)$
factoring through the pair groupoid $P(S^2)$.

Then we have a 2-commutative diagram
\[
 \xymatrix{P_0 TS^2 \ar[d] \ar[rd] \ar[rr]^{\phi_{1}} & & P_0 A \ar[d] \\
\cG(TS^2) \ar[r]^{pr} & P(S^2) \ar[r]^{\Phi_{gpd}} & \cG(A)},
\]where $ \phi_{1}$ is the natural map induced by $\phi_{algd}$
between the path spaces. The triangle 2-commutes by the definition of $pr$. The big square 2-commutes since $\tilde{\Phi}_{gpd}$ has its  infinitesimal morphism  $\phi_{algd}$. Thus the right tetrahedron 2-commutes. 

By Lemma \ref{lemma:strict-hs}, $\Phi_{gpd}$ is presented by a strict
groupoid morphism $\phi_{gpd}: Hol(P_0 TS^2) = P_0 TS^2 \times_{S^2\times S^2}
P_0 TS^2 \xrightarrow{\phi_1 \times_{\Phi_{gpd}} \phi_1} P_0 A
\times_{\cG(A)} P_0 A = Mon(P_0 A)$. 
Thus on the level of Lie groupoids presenting the differentiable stacks of
the arrow space, we have strict maps of Lie groupoids,
\[
\begin{CD}
\tphi_{gpd}: Mon(P_0 TS^2) @>pr>> Hol(P_0 TS^2) @>\phi_{gpd}>> Mon(P_0 A).
\end{CD}
\]

A Lie algebroid morphism $TS^2 \to A$ is an $A$-homotopy such that it
gives a loop in $P_a A$ along the foliation $\cF$ (see Section
\ref{sec:stacky-v}). Thus a re-parametrization\footnote{We can always
  reparametrize an $A$-path so that it becomes an $A_0$-path. See also \eqref{eq:tau}.} of  
the map $id: TS^2\to TS^2$ gives a loop in $P_0 A$ based at the constant path $0_N$, thus representing an element 
$\epsilon \in Mon(P_0 TS^2)$.  However 
$pr(\epsilon)=1_{0_N} \in Hol(P_0
TS^2)$ because any
groupoid Morita equivalent to a manifold must have trivial
isotropy groups.  Therefore
\[\tphi_{gpd}(\epsilon)=\phi_{gpd}(1_{0_N})=1_{\phi_1(0_N)}=1_{0_x}. \]
Notice that $\tphi_{gpd} (\epsilon)\in Mon(P_0 A)$ is represented by a
re-parametrization of the $A$-homotopy $\phi_{algd}: TS^2
\to A$, and $1_{0_x}$ is represented by the constant $A$-homotopy
$TS^2 \xrightarrow{0_x} A$.  A re-parametrization of $\phi_{algd}$ is
connected to $\phi_{algd}$ by a path of $A$-homotopies because a
representation of any element in $S_n(A)$ is homotopic to itself via
an element in $S_{n+1} (A)$ \cite[Remark 3.6]{brahic-zhu}. 
 What we mod
out to form the monodromy groupoid $Mon(P_0 A)$ is paths of
$A$-homotopies.  By \cite[Lemma 4.5]{brahic-zhu}, a path of
$A$-homotopies between $\phi_{algd}$ and the constant $A$-homotopy $0_x$ gives a Lie algebroid morphism
$TD^3 \to A$.  We then
arrive at a higher homotopy $F_\phi: TD^3 \to A$ such
that its restriction to the tangent bundle of the boundary of $D^3$ is $F_\phi|_{TS^2}=\phi_{algd}$.

Then the rest is clear: $F_\phi$ integrates to a W-groupoid morphism $F_\Phi:
P(D^3)\to \cG(A)$. Its restriction $F_\Phi|_{P(S^2)}$
corresponds infinitesimally to the Lie algebroid morphism
$\phi_{algd}$ since $F_\phi|_{TS^2}=\phi_{algd}$. By Theorem \ref{thm:lieii}, we have $F_\Phi|_{P(S^2)}\sim\Phi_{gpd}$.
Moreover $F_\Phi|_{D^3 \times N} $ fixes the base point $N\times N$.  Therefore
$F_\Phi|_{D^3\times N}$ gives a homotopy between $\varphi(\sim
F_{\Phi}|_{S^2 \times N})$ and the
constant map $i_x = F_\Phi|_{N\times N }  $. Hence the source fibre $\bs^{-1}(x)$ is
2-connected.
\end{proof}

\begin{lemma} \label{lemma:strict-hs}
Suppose that a morphism of differentiable stacks $\cX \xrightarrow{\Phi} \cY$ and a morphism of their charts $X_0\xrightarrow{\phi}Y_0$ fit into the following 2-commutative diagram 
\begin{equation}\label{diag:alpha}
\xymatrix{ X_0 \ar[r]^{\phi} \ar[d]_{p_x} & Y_0 \ar@{=>}[dl]_{\alpha} \ar[d]^{p_y}\\
\cX \ar[r]^{\Phi} & \cY,
}
\end{equation} where $p_x: X_0 \to \cX$ and $p_y: Y_0 \to \cY$ are the chart
projection maps. Then the H.S. morphism presenting $\Phi$ between the groupoid $X$ and
$Y$ is a strict Lie groupoid morphism.
\end{lemma}
\begin{proof}
We form the following 2-commutative diagram
\begin{equation}\label{diag:xy}
\xymatrix{& Y_1 \ar[rr] \ar[dd] & & Y_0 \ar[dd]^{p_y} \\
X_1 \ar@{.>}[ur]_{\phi_1} \ar[dd] \ar[rr]& & X_0 \ar[ur]_{\phi}
\ar[dd]^{p_x} \\
& Y_0  \ar[rr]_{p_y}  & & \cY.  \\
X_0 \ar[ur]_{\phi} \ar[rr]_{p_x} & & \cX \ar[ur]^{\Phi} }  \end{equation}

There are two composed maps $X_1 \to X_0 \to Y_0 \to \cY  $, one goes
through the upper $Y_0$ and the other goes through the lower
$Y_0$. These two maps are the same up to a 2-morphism since  the
front, back, right, and bottom
faces of the diagram  are 2-commutative. Thus by the property of
pull-backs, there exists a  morphism $\phi_1: X_1 \to Y_1$ making
all the faces of the diagram 2-commutative. This map can be chosen
canonically as in \eqref{eq:phi-1} after we describe fibre products of
stacks in more detail: Given $U\in \cC$ the category of differential
manifolds,  $X_1=X_0\times_{p_x, \cX, p_x} X_0$ as a stack is the category fibred over $\cC$ with
\begin{description}
 \item[objects over $U$:]$(a, \gamma, a')$, where  $a:U\to X_0, a': U\to X_0$ are objects in $X_0$ viewed as a stack, 
$ \gamma: p_x(a')\to p_x(a)$ is a morphism in $\cX$; 
\item[morphisms over $U\xrightarrow{\phi}V$:]  $(a, \gamma, a') \to (b, \kappa, b')$ which is made up by morphisms  $f: a\to b$ and $f': a'\to b'$ in $X_0$ over $U\xrightarrow{\phi} V$ such that 
\[\xymatrix{p_x(a) \ar[d]_{p_x(f)} &\ar[l]_{\gamma} p_x(a') \ar[d]_{p_x(f')}\\
p_x(b) &\ar[l]_{\kappa} p_x(b')
}
\]commutes.
 \end{description}

Then the map $\phi_1$ can be explicitly described as 
\begin{equation} \label{eq:phi-1}
\begin{split}
&(a, \gamma, a') \mapsto \big(\phi(a), \alpha_a^{-1}\circ \Phi(\gamma) \circ \alpha_{a'}, \phi(a') \big)\\
&(f, f') \mapsto \big(\phi(f), \phi(f') \big), \end{split}
\end{equation}on the level of objects and morphisms respectively. Let us now explain how this formula makes sense. On the object level, $\phi_1$ is well-defined because the 2-morphism $\alpha: p_y \circ \phi \Rightarrow \Phi\circ p_x $ fits in the diagram
\[
 \xymatrix{p_y \circ \phi(a) \ar[d]_{\alpha_a} & p_y\circ \phi(a') \ar[d]_{\alpha_a'} \\
\Phi\circ p_x(a) & \ar[l]_{\Phi(\gamma)} \Phi\circ p_x(a')
}
\]
On the morphism level, $\phi_1$ is well-defined because  the following
diagram commutes by the naturality of $\alpha$: 
\[
 \xymatrix{p_y \phi(a) \ar[d]_{p_y \phi (f)} \ar[r]^{\alpha_a} & \Phi p_x(a) \ar[d]_{\Phi p_x(f)} & \Phi p_x (a') \ar[l]_{\Phi(\gamma)} \ar[d]_{\Phi p_x(f')} & p_y\phi(a') \ar[l]_{\alpha_{a'}} \ar[d]_{p_y\phi (f')} \\
p_y\phi(b) \ar[r]_{\alpha_b} & \Phi p_x (b) & \Phi p_x (b') \ar[l]^{\Phi(\kappa)} & p_y \phi(b') \ar[l]^{\alpha_{b'}} 
}
\]

Now we claim that $\phi_1$ and $\phi$ together form a groupoid morphism. First of all, we recall  the groupoid structure \cite{bx} on $X_1 \rra X_0$. The source and target are simply the right and left projections to $X_0$ respectively; the multiplication on $X_1$ is defined as 
\[ (a, \gamma, a') \cdot (a', \gamma', a'') = (a, \gamma \circ \gamma', a''), \quad (f, f') \cdot (f', f'') = (f, f'') \]
where $\gamma \circ \gamma'$ is simply the composition of morphisms in $\cX$; the identity element is $(a, id_{p_x(a)}, a)$ for each object $a$ in $X_0$. 

Then $\phi_1$ obviously preserves the source and target and
multiplication.  To verify that $\phi_1$ also preserves the identity
element, one uses the fact that functors  preserve identity
morphisms. Then the inverse is automatically preserved. Thus $(\phi_1,
\phi)$ is a groupoid morphism. Then $(\phi_1, \phi)$ induces a stack
morphism $\Phi'$ (unique up to a unique 2-morphism) which makes the
diagram \eqref{diag:alpha} 2-commute. In fact the morphism $\Phi$ in
diagram \eqref{diag:alpha} is totally determined by $\phi$ on the
image of $p_x$, which is a sub-prestack of $\cX$. But since $p_x$ is
an epimorphism, the stackification of the prestack $\im(p_x)$ is
isomorphic to $\cX$ by \cite[52.]{metzler}. Thus any possible $\Phi$
making diagram \eqref{diag:alpha} commute  is unique up to a unique 2-morphism by \cite[51.]{metzler}. Hence the groupoid morphism corresponds to the stack morphism $\Phi$ up to a 2-morphism.

If one prefers to use the language of H.S. morphism, here is a hint at
how 
the H.S. morphism presenting $\Phi$ between $X$ and $Y$ comes from the strict morphism $(\phi_1, \phi)$. The H.S. morphism is
\[X_0 \times_{p_x \circ \Phi, \cY, p_y} Y_0 \cong X_0 \times_{p_y \phi, \cY, p_y} Y_0 \cong X_0 \times_{\phi, Y_0, t} Y_0 \times_{p_y, \cY, p_y} Y_0 \cong X_0 \times_{\phi, Y_0, t} Y_1. \]
\end{proof}

\begin{remark}
One might expect that Theorem \ref{thm:pi2} would give a new proof of
the fact that the simply connected Lie group $G$ has $\pi_2(G)=0$.  By
Theorem \ref{thm:pi2}, we only have $\pi_2(\cG(\g)_1)=0$. Now
$\cG(\g)_1$ is presented by the groupoid $Mon(P_0\g)$ and  the
topological quotient $P_0\g/Mon(P_0\g)_1$ equals $G$. Thus
$\cG(\g)_1=G$ if and only if the monodromy groupoid $Mon(P_0\g)$ has
trivial isotropy groups. This means that for any $a\in P_0\g$, there
is a unique $\g$-homotopy up to higher homotopies between $a$ and
itself. Now $\g$-paths correspond to $G$-paths and $\g$-homotopies
correspond to usual homotopies on $G$. The above statement is
equivalent to the statement that for any $G$-path $g(t)$ there is a unique homotopy up to higher homotopy between $g(t)$ and itself. This exactly requires that $\pi_2(G)=0$. In summary, $\cG(\g)_1=G$ is  equivalent to the fact that $\pi_2(G)=0$. Thus Theorem \ref{thm:pi2} can not give a new proof for $\pi_2(G)=0$. Rather $\pi_2(G)=0$ implies that $G$ is the universal object integrating $\g$ even in the world of 2-groups. This is not true any more for Lie 2-algebras or higher Lie algebras (see \cite[Theorem 7.5]{henriques}).
\end{remark}

\section{Examples}\label{sec:ep}
Now we give some examples of stacky groupoids of the form $\cG(A)$. The purpose of this section is to show that, for
a manifold $M$ with $\pi_2(M)$ nontrivial, $\cG(TM)$ is {\em
not} the traditional homotopy groupoid
$\tilde{M}\times\tilde{M}/\pi_1(M)$ where $\tilde{M}$ is the
simply connected cover of $M$. The second homotopy group $\pi_2(M)$ will play a role too.

\begin{ep} [$\cG(TM)$]          \label{ep:tm} 
The set of Lie algebroid $C^1$-morphisms $TN \to TM$ is equal to the set
of $C^1$-maps $Mor(N, M)$.  With a re-parametrization\footnote{Here to simplify  calculations we will make later, we stop using the 0-boundary condition as in $P_0 A$. We are able to do so because $Mon(P_0 A)$ is Morita equivalent to $Mon(P_a A)$ and both of them present $\cG(A)$ (see \cite[Section 4]{tz:arxiv}).}, the stack $\cG(TM)_1$ is presented by the groupoid $G_1\rra G_0$
with $G_1=Mor(D^2, M)/Mor(D^3, M)$  the space of bigons (the
quotient is similarly given by $\beta_1 \sim \beta_2$ if $\beta_1, \beta_2$ have
the same boundary and bound an element in $Mor(D^3, M)$), and
$G_0=Mor(I, M)$. Here a $k$-dimensional disk $D^k$ is viewed as a $k$-dimensional simplex $\Delta^k$ with many
degenerated faces. If $\pi_2(M)=0$, then two disks $\beta_{1, 2}\in Mor(D^2, M)$ sharing the same boundary bound an element in $Mor(D^3, M)$. If $\pi_1(M)=0$, we see that $G_1=Mor(S^1, M)$. Then two paths $\alpha_{1,2} \in G_0$ sharing the same end points bound exactly one element in $G_1$. Thus $\cG(TM)_1=M\times M$ and $\cG(TM)$ is the pair groupoid $M\times M \rra M$.

\end{ep}

Next we analyse an example where $\pi_2(M)\ne 0$.

\begin{ep}[$\cG(TS^2)$]\label{ep:ts2} We take $M=S^2$ in Example
  \ref{ep:tm}, and we use the same notation as there. We claim that in this case, $G_1 \rra G_0$ is Morita equivalent to the action groupoid $S^3
\times S^3 \times (\R\times \R)/\Z \rra S^3 \times S^3$. Here
$(\R\times \R)/\Z$ is a quotient group by the diagonal $\Z$ action
$(r_1, r_2)\cdot n = (r_1+n, r_2 +n)$, and the action of this
quotient group is given by the projection  $(\R\times \R)/\Z \to
\R/\Z \times \R/\Z =S^1 \times S^1 $ and the product of the $S^1$
action on $S^3$. Our convention of target and source maps are
$\bt(p, q, [r_1, r_2]) = (p,q)$ and $\bs(p, q, [r_1, r_2])=
(p\cdot [-r_1], q \cdot [-r_2])$.

To show this, we give the associated complex line bundle $L\to
S^2$ of the $S^1$-principal bundle $S^3\to S^2$ a Hermitian metric and a compatible connection. We
denote $\xi// \gamma$ as the result of the parallel transport
of a vector $\xi \in L_{\gamma(0)} $ along a path $\gamma$ in $S^2$
to $L_{\gamma(1)}$. Since parallel transport is isometric:
$L_{\gamma(0)} \to L_{\gamma(1)} $, it preserves the $S^1$ bundle
$S^3\subset L$ and the angle $ang(\xi_1, \xi_2)$ between $\xi_1$
and $\xi_2$. Here the angle $ang(\cdot, \cdot)$ $L\oplus L \to
S^1$ is point-wise the usual angular map (or argument map) $\C \to
S^1$. It satisfies
\[ang(\xi_1, \xi_2)+ang(\xi_2, \xi_3) = ang(\xi_1, \xi_3)
\quad \text{and}\quad ang(\xi_1, \xi_2)=-ang(\xi_2, \xi_1). \]
Therefore for two paths $\gamma_1$ and $\gamma_2$ sharing the same
end points, we can define the angle $ang(\gamma_1, \gamma_2)$
between them to be \begin{equation}\label{eq:ang} ang(\gamma_1,
\gamma_2):=ang(\xi//\gamma_1, \xi//\gamma_2), \quad \text{for}\;
\xi \in T_{\gamma_1(0)}S^2.\end{equation} Since parallel
transport preserves the angle, this definition does not
depend on the choice of $\xi$ as the following calculation shows,
\[
\begin{split} ang(\xi_1//\gamma_1, \xi_1//\gamma_2) &=
ang(\xi_1//\gamma_1,
\xi_2//\gamma_1)+ang(\xi_2//\gamma_1,\xi_2//\gamma_2) + ang(\xi_2//\gamma_2,  \xi_1//\gamma_2) \\
&= ang(\xi_1, \xi_2)+ ang(\xi_2//\gamma_1,\xi_2//\gamma_2) +
ang(\xi_2, \xi_1) \\
&= ang(\xi_2//\gamma_1,\xi_2//\gamma_2).
\end{split}
\]
In fact $ang(\gamma_1, \gamma_2)= [\int_{D} \omega_{area}]$ where
$\partial D= -\gamma_1+\gamma_2$ and $\omega_{area}$ is the
standard symplectic (area) form on $S^2$.

\vspace{.6cm}

\begin{equation}\label{pic:s2}
\epsfig{file=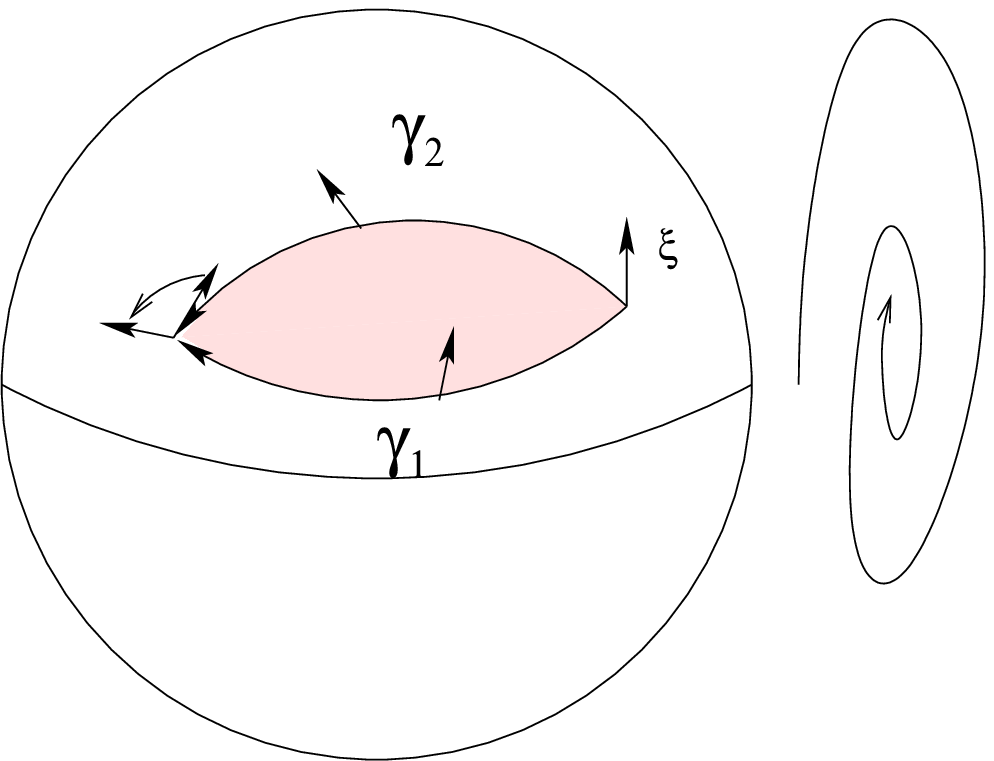,height=3cm}
\end{equation}

As shown in the above picture, because of $\pi_2(S^2)$, $G_1$ is
not simply the pair groupoid $G_0\times_{S^2\times S^2}
G_0$---the set of paths with matched ends. We use
$ang: G_0 \times_{S^2\times S^2}  G_0 \to S^1$ as in
\eqref{eq:ang}, then we have $G_1 = (G_0 \times_{S^2\times S^2}
G_0)\times_{ang, S^1, pr} \R$ where $pr: \R \to S^1$ is the
projection. For example \eqref{pic:s2} corresponds to $(\gamma_1,
\gamma_2, r)$ with $\lfloor r \rfloor = 2$. The pull-back groupoid
of $G_1 \rra G_0$ along the projection (a surjective submersion)
$S^3\times_{\pi,S^2, \bt} G_0 \times_{\bs,S^2, \pi} S^3 \to G_0$
with $\bs(\gamma)=\gamma(0)$, $\bt(\gamma)=\gamma(1)$, and $\pi: S^3\to S^2$,  is
\[ (S^3\times_{S^2} G_0 \times_{S^2} S^3 \times S^3\times_{S^2} G_0 \times_{S^2} S^3 )
  \times_{ang, S^1, pr} \R \rra S^3\times_{S^2} G_0 \times_{S^2}
  S^3. \]

Now  the pull-back groupoid of $S^3 \times S^3 \times (\R\times
\R)/\Z \rra S^3 \times S^3$ along the projection (also a surjective submersion) $S^3\times_{S^2} G_0
\times_{S^2} S^3 \to S^3 \times S^3$ defined by $(\xi,
\gamma,\xi') \mapsto (\xi'//\gamma, \xi//\gamma^{-1})$ is
\begin{equation}\label{eq:pull-back-gpd}
(S^3\times_{S^2} G_0 \times_{S^2} S^3 \times S^3\times_{S^2} G_0
\times_{S^2} S^3 )\times_{(S^3\times S^3) \times (S^3\times S^3) }
S^3\times S^3 \times (\R\times\R)/\Z \rra S^3\times_{S^2} G_0
\times_{S^2} S^3.
\end{equation}
They are the isomorphic as Lie groupoids with the morphism from
the second to the first given by
\[
(\xi, \gamma, \xi') \mapsto (\xi, \gamma', \xi')
\]
on the base of the groupoid, where $\gamma'$ is a path\footnote{Such
  $\gamma'$ can be chosen canonically in the following way. Take the
  usual Riemannian metric on $S^2$. Let $\phi_t(s)$ be the trajectory
  of $\dot{\gamma}(t)^{\perp}$. Then when we increase the angle
  $ang(\dot{\gamma}'(t), \dot{\gamma}(t)//\phi_t(s))$ starting from 0,
  there is a unique $\gamma'$ such that $ang(\gamma, \gamma')=\alpha$
  for any $\alpha \in S^1$. } sharing the same end points of $\gamma$
and such that $ang(\gamma, \gamma') = ang (\xi, \xi'//\gamma)$.
On the level of morphisms, we define
\begin{equation}\label{eq:arrow}
(\xi_1, \gamma_1, \xi'_1, \xi_2, \gamma_2, \xi'_2, [(r_1, r_2)]) \\
\mapsto (\xi_1, \gamma'_1, \xi'_1, \xi_2, \gamma'_2, \xi'_2,
r_2-r_1).
\end{equation} We need to verify that this map does land on the correct manifold.
By \eqref{eq:pull-back-gpd} we have $(\xi_1//\gamma_1^{-1})\cdot
[-r_2]=\xi_2//\gamma_2^{-1}$, thus
\[
\begin{split}
 -[r_2]&=ang(\xi_1//\gamma_1^{-1}, \xi_2//\gamma_2^{-1}) \\
&=ang(\xi_1//\gamma_1^{-1}, \xi_2//\gamma_1^{-1}) +
ang(\xi_2//\gamma_1^{-1}, \xi_2//\gamma_2^{-1} ) \\
&= ang(\xi_1, \xi_2) + ang(\gamma_1^{-1}, \gamma_2^{-1}),
\end{split}
\]
hence $ang(\gamma_1, \gamma_2) = [r_2]+ang(\xi_1,\xi_2)$ and
similarly $ang(\gamma_1,\gamma_2) = -[r_1]-ang(\xi'_1,\xi'_2)$.
Therefore
\[ [r_2-r_1]= 2ang(\gamma_1, \gamma_2) -
ang(\xi_1,\xi_2) +  ang(\xi'_1,\xi'_2). \] Since
\[ ang (\xi_1, \xi'_1 //\gamma_1) + ang(\xi'_1 //\gamma_1, \xi'_2
//\gamma_1) + ang(\xi'_2//\gamma_1, \xi'_2//\gamma_2) +ang(\xi'_2
//\gamma_2, \xi_2) = ang(\xi_1, \xi_2), \]
we have
\[- ang(\xi_1,\xi_2) + ang(\xi'_1,\xi'_2) = -ang(\xi_1, \xi'_1//\gamma_1) +
ang(\xi_2, \xi'_2//\gamma_2) -ang(\gamma_1, \gamma_2). \] Hence
\[
\begin{split}
[r_2-r_1]&=ang(\gamma_1, \gamma_2) + (-ang(\xi_1,
\xi'_1//\gamma_1) + ang(\xi_2, \xi'_2//\gamma_2))\\
&=ang(\gamma_1, \gamma_2) +ang(\gamma'_1, \gamma_1) +
ang(\gamma_2, \gamma'_2) \\
&=ang(\gamma'_1, \gamma'_2)
\end{split}
\]
Therefore \eqref{eq:arrow} is well-defined. It is not hard to
see that it is indeed a groupoid isomorphism.

Therefore $G_1\rra G_0$ and $ S^3\times S^3 \times (\R\times
\R)/\Z \rra S^3 \times S^3$ are Morita equivalent via this third
groupoid \eqref{eq:pull-back-gpd}. Therefore $\cG(TS^2)$ is presented by $ S^3 \times S^3
\times (\R\times \R)/\Z \rra S^3 \times S^3$.
\end{ep}

\begin{remark}[on 2-connected covers] \label{rk:2-conn-cover}
It is well-known that the cohomology group $H^2(M, \Z)$ corresponds to
$S^1$-bundles via Chern classes. Let $B\Z$ denote the differentiable stack corresponding to the Lie group $\Z$. Then $B\Z$ is a stacky Lie group with the multiplication given by the group(oid) morphism $\Z \times \Z \xrightarrow{+} \Z$ on the level of presenting groupoids. The morphism $+$ is a groupoid morphism because $\Z$ is abelian.  Then from a stacky viewpoint, $H^2(M, \Z)$
also corresponds to $B\Z$-gerbes: Given a covering $U_i$ of $M$, and  $g_{ijk}$ representing a
\v{C}ech class in $H^2(M, \Z)$,  there is a stack $\cG$ presented by groupoid
$\sqcup U_{ij} \times \Z \rra \sqcup U_i $ with the groupoid
multiplication $(x_{ij}, n)\cdot (x_{jk}, m) = (x_{ik},
n+m+g_{ijk})$ and the source and target maps inherited from the
groupoid $\sqcup U_{ij} \rra \sqcup U_i $ which presents $M$. Here
$x_{ij}$ denotes a point in $U_{ij}$.
An $A$-gerbe corresponds to an $A$-groupoid central extension for
an abelian group $A$ \cite{bx}.  It is clear that $\cG$  is a $\Z$-gerbe (or
equivalently a $B\Z$-principal bundle \cite{wz}) over $M$ from the
following diagram of groupoid central extension:
\[ \xymatrix{ 1 \ar[r] & \Z \times \sqcup U_i \ar[r] \ar[d]\ar@<-1ex>[d] & \sqcup U_{ij} \times
\Z \ar[r] \ar[d]\ar@<-1ex>[d]  & \sqcup U_{ij} \ar[r]
\ar[d]\ar@<-1ex>[d] & 1 \\
&  \sqcup U_i \ar@{=}[r]  & \sqcup U_i \ar@2{-}[r]  & \sqcup U_i &
}\] 

For example,  corresponding to the Hopf fibration $S^3\to S^2$, there is
a $\Z$ gerbe or $B\Z$ principal bundle denoted as $\tS^2$. Now,  apply
the long exact sequence of homotopy groups to the $B\Z$-fibration
$\tS^2 \to S^2$. Since $\pi_{n\ne 1}(B\Z)=0$ and $\pi_1(B\Z)=\Z$, we
have $\pi_1(\tS^2)=\pi_2(\tS^2)=0$ and $\pi_{\geq 3} (\tS^2) =
\pi_{\geq 3} (S^2)$. Hence we can view $\tS^2$ as a 2-connected
``covering'' of $S^2$. Although $S^3$ also has the same property of
homotopy groups, $\tS^2$ seems to be a more suitable candidate than
$S^3$ since $\tS^2$ has the same dimension as $S^2$ given $\dim
B\Z=0$. 

In this  example, $S^3 \times \R \rra S^3$ presents the stack
$\tilde{S^2}$.  What we prove above shows that the stacky groupoid
$\cG(TS^2)$ is not the usual one, $S^2 \times S^2$, but rather $\tilde{S^2}
\times \tilde{S^2}/B\Z$. This is to ensure
that it has 2-connected source-fibre, in agreement with the  general result we
proved in Theorem \ref{thm:pi2}. This resembles the construction of the
groupoid integrating $TM$ when $M$ is a non-simply-connected
manifold. The source simply connected groupoid of $TM$ is
$\tilde{M}\times \tilde{M}/\pi_1(M)$.  Here, for $\cG(TS^2)$, the
construction is comparable to this, but on a higher level---the
aim is to kill $\pi_2$ of the source fibre. Further, using the long exact
sequence of homotopy groups, we have $\pi_2(\tS^2\times
\tS^2/B\Z)=\pi_1(B\Z)=\Z$ since $\pi_2(\tS^2\times \tS^2)=0$.  
\end{remark}

\bibliographystyle{habbrv}
\bibliography{../bib/bibz}

\end{document}